\documentclass[11pt]{article}

\ifx\comments\undefined
  \usepackage[letterpaper,margin=1in,includefoot]{geometry}
  \newcommand{\XXXcomment}[1]{}
  \newcommand{\XXXcommentR}[1]{}
  \newcommand{\XXXcommentG}[1]{}
\else
  \usepackage[letterpaper,margin=.5in,right=2.5in,marginpar=2in,includefoot]{geometry}
  \newcommand{\XXXcomment}[1]{\marginpar{\color{blue}{\footnotesize #1}}}
  \newcommand{\XXXcommentR}[1]{\marginpar{\color{red}{\footnotesize #1}}}
  \newcommand{\XXXcommentG}[1]{\marginpar{\color{green}{\footnotesize #1}}}
\fi

\ifx\arXiv\undefined
  \usepackage[pdftex,
              bookmarks=true,
              bookmarksnumbered=true,
              pdfborder={0 0 0},
              plainpages=false,
              unicode,
              pdfpagelabels]{hyperref}
\fi

\usepackage{amssymb,amsthm,latexsym,amsmath,mathrsfs,stmaryrd,graphicx,fancyhdr,paralist,longtable, multirow, wasysym}

\usepackage[x11names]{xcolor}
\usepackage{colortbl}

\usepackage[compact]{titlesec}
\usepackage[titletoc,title]{appendix}
\usepackage{multicol}

\usepackage{bussproofs}

\usepackage{linguex}

\ifx\arXiv\undefined
  \usepackage[sort]{natbib}
\fi

\usepackage{caption,float}
\usepackage{tikz}
\usetikzlibrary{trees,arrows,automata,shapes,decorations,decorations.pathmorphing}

\usepackage[normalem]{ulem}

\usepackage[makeroom]{cancel}

\theoremstyle{definition} 
\newtheorem{definition}{Definition}[section]
\newtheorem{theorem}[definition]{Theorem}

\newtheorem{conjecture}[definition]{Conjecture}
\newtheorem{corollary}[definition]{Corollary}
\newtheorem{example}[definition]{Example}

\newtheorem{fact}[definition]{Fact}
\newtheorem{counterexample}[definition]{Counter-example}

\definecolor{darkgray}{rgb}{.1, .1, .14}
\newenvironment{myproof}{\footnotesize\begin{proof}\color{darkgray}}{\end{proof}\normalsize}

\newtheorem{thm}{Theorem}

\newcommand{\ourtitle}{Extracting Conditionals from Consequence Relations}

\newcommand{\paul}{Paul Egr\'e}

\newcommand{\paulinstitute}{IJN / ENS}

\newcommand{\emmanuel}{Emmanuel Chemla}

\newcommand{\emmanuelinstitute}{LSCP / ENS}

\newcommand{\emmanuelFunding}{Funded by the European Research Council under the European Union's Seventh Framework Programme (FP/2007-2013) / ERC Grant Agreement n. 313610 and was supported by ANR-10-IDEX-0001-02 PSL* and ANR-10-LABX-0087 IEC.}
  
\newcommand{\paulFunding}{Funded by Programs ANR-10-LABX-0087 IEC,
  ANR-10-IDEX-0001-02 PSL*,  by the ANR project ``Trivalence and
  Natural Language Meaning'' (ANR-14-CE30-0010), and by Ministerio de Econom\'ia, Industria y Competitividad, Gobierno de Espana as part of the project ``Logic and substructurality", Grant. n FFI2017-84805-P.}

\ifx\arXiv\undefined
\hypersetup{ 
  pdfauthor={\emmanuel{}, \paul{}}, 
  pdftitle={\ourtitle{}}
}
\fi

\author{
  \emmanuel{}\footnote{\emmanuelFunding{}}\\
  \small\emmanuelinstitute{}
  \and
  \paul{}\footnote{\paulFunding{}}\\
  \small\paulinstitute{}
 }

\date{}

\def\truthrelation{\mathrel|\joinrel\equiv}

\usepackage[T1]{fontenc}
\usepackage[sc,osf]{mathpazo}

\usepackage{microtype}
\usepackage{nicefrac}

\newcommand{\dint}[1]{{\left\vert\kern-0.25ex\left\vert\kern-0.25ex\left\vert #1 
    \right\vert\kern-0.25ex\right\vert\kern-0.25ex\right\vert}}

    \newcommand{\half}{\nicefrac{1}{2}}

\begin{document}

\renewcommand*{\thefootnote}{\fnsymbol{footnote}}
\noindent {\Large From Many-Valued Consequence to Many-Valued Connectives%
\footnote{Acknowledgements: We are very grateful to Benjamin Spector for providing inspiration and support to this project. We thank two anonymous referees for detailed and helpful comments. We also thank Denis Bonnay, Keny Chatain, Christian Ferm\"uller, Jo\~ao Marcos, Hitoshi Omori, Francesco Paoli, David Ripley, Lorenzo Rossi, Hans Rott, Philippe Schlenker, Jan Sprenger, Shane Steinert-Threlkeld, Heinrich Wansing for helpful conversations, as well as audiences in Regensburg (workshop ``New Perspectives on Conditionals and Reasoning''), Bochum  (Logic in Bochum IV) and Dagstuhl (Dagstuhl Seminar 19032 ``Conditional logics and conditional reasoning''). The research leading to these results has received funding from the European Research Council under the European Union's Seventh Framework Programme (FP/2007-2013) / ERC Grant Agreement n.313610, and from the ANR program ``Trivalence and Natural Language Meaning'' (ANR-14-CE30-0010). We also thank the Ministerio de Econom\'ia, Industria y Competitividad, Gobierno de Espana, as part of the project ``Logic and substructurality" (Grant. FFI2017-84805-P), as well as grant FrontCog, ANR-17-EURE-0017 for research conducted in the Department of Cognitive Studies at ENS.
We dedicate this paper to the memory of Carolina Blasio.
}}\\
\renewcommand*{\thefootnote}{\arabic{footnote}}
\setcounter{footnote}{0}\\
Emmanuel Chemla$^\textrm{a}$ \& Paul Egr\'e$^\textrm{b}$ \\
{\scriptsize a. Laboratoire de Sciences Cognitives et Psycholinguistique, D\'epartement d'\'etudes cognitives, ENS, PSL University, EHESS, CNRS, 75005 Paris, France}\\
{\scriptsize b. Institut Jean Nicod, D\'epartement d'\'etudes cognitives \& D\'epartement de philosophie, ENS, PSL University, EHESS, CNRS, 75005 Paris, France}
%Emmanuel Chemla$^\textrm{a,c}$ \& Paul Egr\'e$^\textrm{b,c}$ \\
%{\scriptsize a. LSCP, D\'epartement d'\'etudes cognitives}\\
%{\scriptsize b. Institut Jean Nicod, D\'epartement d'\'etudes cognitives \& D\'epartement de philosophie}\\
%{\scriptsize c. ENS, EHESS, CNRS, PSL University, 75005 Paris, France}

\begin{abstract}
\noindent Given a consequence relation in many-valued logic, what connectives can be defined? For instance, does there always exist a conditional operator internalizing the consequence relation, and which form should it take? In this paper, we pose this question in a multi-premise multi-conclusion setting for the class of so-called intersective mixed consequence relations, which extends the class of Tarskian relations. Using computer-aided methods, we answer extensively for 3-valued and 4-valued logics, focusing not only on conditional operators, but also on what we call Gentzen-regular connectives (including negation, conjunction, and disjunction). For arbitrary $N$-valued logics, we state necessary and sufficient conditions for the existence of such connectives in a multi-premise multi-conclusion setting. 
The results show that mixed consequence relations admit all classical connectives, and among them pure consequence relations are those that admit no other Gentzen-regular connectives. Conditionals can also be found for a broader class of intersective mixed consequence relations, but with the exclusion of order-theoretic consequence relations.
\end{abstract}

\noindent{\small  {\bf Keywords:} logical consequence; mixed consequence; compositionality; {truth-functionality}; many-valued logic; algebraic logic; conditionals; connectives; sequent calculus; deduction theorem; truth value}\bigskip

\section{Introduction: matching conditionals and consequence relations}

In 2-valued logic, logical consequence is defined as the preservation of the value True (=$1$) from premises to conclusions in an argument. A fundamental feature of 2-valued logic is the existence of a binary sentential connective {representing or internalizing} the consequence relation in the object-language, namely the material conditional $\supset$ (or material implication, taking the value 1 exactly when the value of the antecedent is less or equal to the value of the consequent). As established by the \emph{deduction theorem}, $A\vdash B$ iff {$A\supset B$ is valid}, and more generally, when $\Gamma$ is a finite set of premises $A_{1}$ to $A_{n}$: $\Gamma\vdash B$ iff {$(A_{1} \supset (A_{2} \supset \cdots (A_{n} \supset B))\cdots)$ is valid}.
{This work extends the question of the representability of consequence relations by means of adequate conditional operators to logics with more truth-values, following a question originally raised by Benjamin Spector and posed in our joint article \cite{chemla2017charac}.
}

In many-valued logics, more notions of consequence become available as the set of truth values expands, and similarly the space of binary operators quickly increases. Perhaps as a result of the greater freedom in the choice of those parameters -- though sometimes due to {independent} desiderata -- various popular systems of many-valued logics rest on the choice of a conditional operator that fails one or the other direction of the deduction theorem, and thereby fails to internalize logical consequence adequately in the object-language (see \citealt{avron1991natural, cobreros2015vagueness, wintein2016kleene}). 
{For instance, a popular logic such as Strong Kleene's $\mathsf{K3}$, in which $v(A\rightarrow B)=v(\neg A \vee B)= max(1-v(A), v(B))$, and where consequence is defined as the preservation of the value 1, loses conditional introduction ($A\vdash A$, but $\not\vdash A\rightarrow A$). 
Similarly, the dual logic $\mathsf{LP}$ which uses the same conditional operator but in which consequence is defined as the preservation of non-zero values ($\{1, \half\}$) fails the converse direction akin to the modus ponens ($\vdash A \rightarrow ((A \rightarrow B) \rightarrow B)$, but $ A, (A\rightarrow B)\not\vdash B$). {Both of those counterexamples are forestalled in \L ukasiewicz's 3-valued system $\mathsf{\L3}$, who defines logical consequence as the preservation of the value $1$, and in which $v(A\rightarrow B)$ equals $1$ whenever $v(A)\leq v(B)$, and equals $1-(v(A)-v(B))$ otherwise (\citealt{luk1920}). Even there, however, conditional introduction fails, for whereas $A\wedge \neg A\vdash B$, the schema $A\wedge \neg A\to B$ is not valid (as can be seen by assigning $A$ the value $\half$ and $B$ the value $0$).}\footnote{{See \cite{pogorzelski1964deduction} and \cite{avron1991natural} for more on the deduction theorem in relation to \L ukasiewicz's three-valued conditional. Pogorzelski shows that \L ukasiewicz's conditional satisfies a more complex form of the deduction theorem relative to consequence defined as the preservation of the value 1. Avron shows that it can satisfy the deduction theorem in standard form if the definition of logical consequence is modified in a way that rules out the above counterexample, by giving up structural contraction. We note that relative to the mixed consequence relation $st$ (introduced below), \L ukasiewicz's conditional would satisfy the deduction theorem (but not what we call premise Gentzen-regularity, unlike with Avron's consequence)}.}
}

This situation of mismatch between consequence and conditional raises two natural questions: (i)~Given a consequence relation in $N$-valued logic, which conditional operators satisfy the deduction theorem, if any? The question, we will see, does not always have an obvious answer for $N \geq 3$. (ii)~Moving away from conditional operators, which connectives in general are naturally associated with a consequence relation, and what do the corresponding operators say, taken together, about this consequence relation? The first issue provides the motivation and a guiding direction for this paper, but we provide important general results pertaining to question~(ii). {By so doing, we follow the general perspective that \cite{bonnay2012consequence} aptly call ``extracting or mining'' connectives from consequence relations. 
% One difference between the projects may be described by noting that we start from the consequence relation defined intensionally, while they start from a consequence relation defined extensionally over a language

As a first step toward that goal, we are interested in a wider space of consequence relations than those standardly defined in terms of the preservation of a constant set of designated values from premises to conclusion. One motivation for this comes from the consideration of so-called \emph{mixed} consequence relations (\citealt{cobreros2012tolerant, chemla2017charac, frankowski2004, malinowski1990q}). By way of illustration, consider the semantics of the Strong Kleene conditional given above: re.~(ii), we just observed that it internalizes logical consequence neither as the preservation of the value $1$ (so-called $ss$-consequence), nor as the preservation of non-zero values (so-called $tt$-consequence). However, it does adequately internalize so-called \emph{strict-tolerant} consequence ($st$), defined in terms of the impossibility for the premises to take the value $1$ and for the conclusion to take the value $0$ (see \citealt{cobreros2012tolerant, cobreros2015vagueness}).

In what follows, we thus propose to investigate the relationship between consequence relations and conditional operators over the space of \emph{intersections of mixed consequence relations}. As established in previous work (see \citealt{chemla2018suszko}), intersective mixed relations fundamentally correspond to monotone consequence relations (see also \citealt{blasio2017inferentially,french2017valuations}). This space includes non-Tarskian relations which may fail either to be reflexive ($ts$) or transitive ($st$). On the other hand, it does include some Tarskian relations that are not obtainable in terms of the preservation of a fixed set of designated values (see \citealt{chemla2017charac}). This includes in particular the \emph{order-theoretic} consequence relations (definable as the intersection of all pure relations), an example of which is $ss\cap tt$ in 3-valued logic (requiring the preservation of both 1 and of non-zero values from premises to conclusion). While previous work exists concerning the relationship between Tarksian consequence operations and the existence of conditional operators internalizing them (see \citealt{avron1991natural}), we are not aware of a similar systematic investigation over the wider space of relations here considered.\footnote{A recent exception is \cite{wintein2016kleene} looking at 3-valued and 4-valued mixed consequence relations, but not at intersective mixed relations.} 

Before we start, let us make a few more remarks concerning our research agenda. Firstly, like \cite{avron1991natural}, in this paper we will consider not just conditional operators, but more logical connectives, in particular negation, conjunction and disjunction, and we will be interested in the internalization of consequence relations by operators that not only satisfy the deduction theorem, but that satisfy further conditions that seem just as natural.

Secondly, one of the ambitions of this paper is to serve as a repository of results that can help logicians to find out a map of consequence relations and conditionals for the special case of 3-valued and 4-valued logics in particular. There is a sense in which every logician knows the map of 2-valued logic: they are completely familiar with all 16 binary operators, and they know of several arguments to select the horseshoe as the best candidate for being a conditional in that space. For 3-valued logic, where the number of binary truth-functional operators approximates 20,000, no similar map is available, and even for the best-known logics, it can be unclear which set of connectives is to be paired with a given consequence relation, or conversely. We propose to fill this gap, namely to chart the land of 3-valued and 4-valued logic in a systematic manner. %By so doing, we follow the general perspective that \cite{bonnay2012consequence} aptly call ``extracting or mining'' connectives from consequence relations. 

Thirdly, to serve that goal we will present several results based on computer-aided methods. This means that for several of the results we will state regarding 3-valued and 4-valued logic, we have found useful to do an exhaustive search of the space of consequence relations and the associated binary operators. 
The search was used both as a heuristic to discover generalizations as well as a way to demonstrate results, and we think such computer-aided exploration can be of value to answer further questions.\footnote{Computer-aided investigations of this kind still seem quite rare, which is striking considering that some pioneers such as \cite{Foxley1962:computer3valued} had bravely started deploying them for very related tasks, when much more ingenuity was needed to compensate for the lower power of computers.}

The paper is structured as follows. 
In Section~\ref{sec:consequence}, we lay out the ground for the rest of the paper: we define consequence relations and introduce \emph{intersective mixed consequence truth-relations} as our framework. 
In Section~\ref{sec:operators}, we put forward the notion of a regular connective, and focus on what we call \emph{Gentzen-regular connectives}, that is connectives obeying a biconditional version of Gentzen's classic sequent calculus rules. We show that classical connectives are Gentzen-regular, but that Gentzen-regular connectives are more general: they correspond to {a subset of} all truth-functional connectives in a given 4-valued logic. The problem we pose in this paper, generally put, comes down to determining which class of Gentzen-regular connectives is admitted by a given intersective mixed consequence truth-relation. 
In Sections~\ref{sec:three} and~\ref{sec:four}, we list extensively which consequence relations in 3-valued and 4-valued logics admit conjunctions, disjunctions, negations and conditionals, providing computer programs to reproduce and extend this inventory. 
In Section~\ref{sec:pureandorderresults}, we move on to $N$-valued logics and show that the problem admits a simple and stable solution for two important classes of consequence relations, \emph{mixed} consequence relations (and even more so for the subclass of \emph{pure} consequence relations) and \emph{order-theoretic} consequence relations. Finally, we lay out algebraic characterization results for all finite-valued logics in Section~\ref{sec:N}.

\section{Consequence relations}\label{sec:consequence}

In this section, we introduce the framework we will be using to represent and interpret consequence relations and connectives. 

\subsection{Languages and semantics}

We work with sentential languages, and restrict attention to truth-functional interpretations of formulae. Throughout the paper, we use Roman capitals $A, B, ...$ to denote formulae of the language, and Greek capitals $\Gamma$, $\Delta$,..., to denote sets of formulae. We use small Greek letters, $\gamma, \delta$, to denote subsets of truth values.

\begin{definition}[Sentential Language] A \emph{sentential language} $\mathcal{L}$ consists of a denumerable set of atoms $p_{1}, p_{2},...$, together with a set $\mathcal{C}$ of sentential connectives, where formulae are generated in the usual way. 
\end{definition}

\begin{definition}[Semantic interpretation]

Given a sentential language $\mathcal{L}$, and a set of truth values $\mathcal{V}$ containing at least the special values $1$ and $0$, a \emph{semantic interpretation} (or valuation) is a morphism $v$ such that
\begin{itemize}
\item $v$ maps atoms on truth values 
\item $v$ maps $n$-ary sentential connectives (a.k.a. operators) to truth-functions from $\mathcal{V}^{n}$ to $\mathcal{V}$
\item for any $n$-ary connective $C$, $v(C(A_{1},...,A_{n}))=v(C)(v(A_{1}),...,v(A_{n}))$.
\end{itemize}
\end{definition}

\begin{definition}[Semantics]\label{def:semantics}

A \emph{semantics} is a set of valuations such that
\begin{itemize}
\item for any connective $C$ and any two valuations $v_1$, $v_2$: $v_1(C)=v_2(C)$,
\item for any list of pairs of (distinct) atoms and truth values, one can find a valuation $v$ which assigns to each of the atoms in the list the relevant truth value.
\end{itemize}

\end{definition}

In this paper, unless otherwise noted, we will also suppose that the semantics is `sufficiently' expressive. In the following Definition~\ref{def:constantsandmaximalexpressiveness}, we present three levels of expressiveness: \emph{maximal expressiveness}, \emph{constant expressiveness}, and \emph{atomic expressiveness}. The first one is the most stringent, the second one may be obtained by ensuring the presence of relevant $0$-ary connectives in the language. The last one, atomic expressiveness, is the least demanding. We call it atomic expressiveness because it holds for instance if the language has atomic formulae whose semantic values are meant to cover the whole space of truth-values and with maximal variation across the different atomic formulae, as in the columns of a truth-table (this situation corresponds to a `valuational' semantics in the sense of \citealp{chemla2018suszko}).
Maximal expressiveness and constant expressiveness play a useful role in the rest of the paper, and by default we assume that they hold. We will make it explicit when results hold for the more inclusive class of atomic expressive semantics.

\begin{definition}[Maximal expressiveness, Constant expressiveness, Atomic expressiveness]\label{def:constantsandmaximalexpressiveness}
$\phantom{x}$
\begin{itemize}
\item The semantics is \emph{maximally expressive} if for every function $P$ from the set of valuations to the set of truth values (that is, for every `proposition'), the language contains a formula $F_P$ such that for all valuation $v$: $v(F_P)=P(v)$.
\item The semantics is \emph{constant expressive} if for every truth value $\alpha$, the language contains a formula $F_\alpha$ such that for all valuation $v$, $v(F_\alpha)=\alpha$. 

\item The semantics is \emph{atomic expressive} if for every set of truth values $\gamma$, there is a set of formulae $\Gamma$ and a valuation $v$ such that $v(\Gamma)=\gamma$.

\end{itemize}
\end{definition}

\subsection{Consequence relations: intersective mixed}

We define a consequence relation to be any relation between subsets of formulae of the language (so potentially a relation between several premises and several conclusions, following \citealt{gentzen1935investigations, scott1974completeness, shoesmith1978multiple}). At the semantic level, we restrict attention to consequence relations that are \emph{truth-relational}, that is, consequence relations interpretable as relations between sets of truth values:

\begin{definition}[Consequence Relation]
Given $\mathcal{L}$ a sentential language,
we call a \emph{consequence relation} a subset $\vdash$ of $\mathcal{P}(\mathcal{L})\times\mathcal{P}(\mathcal{L})$.\end{definition}

\noindent We note that this way of defining consequence relations, through sets rather than lists, imposes structural contraction (namely $\Gamma, A, A\vdash \Delta$ iff $\Gamma, A\vdash \Delta$, and similarly $\Gamma\vdash A, A,\Delta$ iff $\Gamma \vdash A, \Delta$). 

\begin{definition}[Truth-relations, truth-relational consequence relations] A \emph{truth-relation} is a subset of $\mathcal{P}(\mathcal{V})\times\mathcal{P}(\mathcal{V})$. We say that a consequence relation $\vdash$ is \emph{truth-relational} if there exists a truth-relation  $\truthrelation$ such that $\Gamma \vdash \Delta$ iff for every semantic interpretation $v$ in the semantics, $v(\Gamma) \truthrelation v(\Delta)$.
\end{definition}

We will restrict attention to 
specific types of truth-relations, namely \emph{intersective mixed consequence relations}. The notion of a \emph{mixed consequence relation} constitutes a generalization of the classic semantic notion of consequence relation in the sense of Tarski (see \citealt{cobreros2012tolerant, chemla2017charac}). The classic notion is defined as the preservation of designated values from premises to conclusion in an argument. For a mixed consequence relation, the set of designated values is allowed to vary between premises and conclusions:

\begin{definition}[Designated values]
A set of designated values is a set of truth values including $1$ and not $0$. 
\end{definition} 

\begin{definition}[Mixed consequence, pure consequence and intersective mixed consequence relations]\
\begin{itemize}
\item A \emph{mixed consequence truth-relation} is a truth-relation noted $\truthrelation_{\mathcal{D}_p,\mathcal{D}_c}$, where
	$\mathcal{D}_p$ is a \emph{premise-set of designated values}, 
	and $\mathcal{D}_c$ is a \emph{conclusion-set of designated values},  such that 
	for all sets of truth values 
			$\gamma, \delta: \gamma \truthrelation_{\mathcal{D}_p,\mathcal{D}_c} \delta$ 
			iff 
			$\gamma\subseteq\mathcal{D}_p$ implies $\delta\cap\mathcal{D}_c\not=\emptyset$.
\item If the sets of designated values are the same for premise and conclusion, that is $\mathcal{D}_p=\mathcal{D}_c$, the relation is called a \emph{pure consequence relation}.
\item An \emph{intersective mixed consequence truth-relation} is an intersection of mixed consequence relations: $\truthrelation_{\mathcal{D}_p^1,\mathcal{D}_c^1}\cap...\cap\truthrelation_{\mathcal{D}_p^K,\mathcal{D}_c^K}$. 
\end{itemize}
\end{definition}

\begin{definition}[Representation, Minimal Representation]
A \emph{representation} of an intersective mixed consequence truth-relation is a set of mixed consequence relations {whose intersection} is the relation in question.
A \emph{minimal representation} of an intersective mixed consequence truth-relation is a representation based on the least possible number of mixed truth-relations whose intersection gives that relation.
\end{definition}

\begin{example}[Mixed relations]

Let $\mathcal{V}=\{0, \half,1\}$. Let $\mathcal{D}_{s}=\{1\}$, and $\mathcal{D}_{t}=\{1,\half\}$. Then $\truthrelation_{\mathcal{D}_{s}, \mathcal{D}_{s}}$, also called $ss$, and $\truthrelation_{\mathcal{D}_{t}, \mathcal{D}_{t}}$, also called $tt$, are mixed consequence relations corresponding to standard consequence relations (also called \emph{pure} relations, because the premise and conclusion sets of designated values are identical, see \citealt{chemla2017charac}). The relation $\truthrelation_{\mathcal{D}_{s}, \mathcal{D}_{t}}$ is a mixed consequence relation, also known as {p-consequence} (\citealt{frankowski2004}), or $st$ (\citealt{cobreros2012tolerant}), and likewise for $\truthrelation_{\mathcal{D}_{t}, \mathcal{D}_{s}}$, {also known as q-consequence} (\citealt{malinowski1990q}), or $ts$ (\citealt{cobreros2012tolerant}). 
\end{example}

\begin{example} [Intersective Mixed relations]\label{ex:sstt}

Consider the intersective mixed relation given by $\truthrelation_{\mathcal{D}_{s}, \mathcal{D}_{s}}\cap \truthrelation_{\mathcal{D}_{t}, \mathcal{D}_{t}}$. This relation, also known as $ss\cap tt$, cannot be expressed as a mixed relation (for a proof, see \citealt{chemla2017charac}). The latter representation is a minimal representation for it. 
A nonminimal representation for it is for example: $\truthrelation_{\mathcal{D}_s, \mathcal{D}_s}\cap \truthrelation_{\mathcal{D}_t, \mathcal{D}_t} \cap \truthrelation_{\mathcal{D}_s, \mathcal{D}_t}$, i.e. $ss\cap tt \cap st$. (It is a representation for it because $ss\subseteq st$, so $st$ plays no role in the representation, except for making it nonminimal).

\end{example}

Why focus on intersective mixed relations? For three main reasons. First, even aside from our purpose in this paper, mixed consequence relations have had a wide range of applications in recent years (see in particular \citealt{cobreros2015vagueness} for a review). Secondly, as shown in \citealt{chemla2018suszko}, intersective mixed consequence relations correspond exactly to the class of \emph{monotonic} consequence relations, that is, of relations $\vdash$ such that $\Gamma'\vdash \Delta'$ whenever $\Gamma \vdash \Delta$ and $\Gamma\subseteq \Gamma', \Delta\subseteq \Delta'$. This class, importantly, includes consequence relations that are not necessarily reflexive (like $ts$, in which $A\nvdash A$) or transitive (like $st$, in which $A\vdash B$, and $B\vdash C$ need not imply $A\vdash C$). That is, it includes non-Tarskian consequence relations. But it thereby gives a more general perspective on Tarskian relations, by setting reflexivity and transitivity as distinct parameters (on the same perspective, see \citealt{blasio2017inferentially, french2017valuations}). Our third and specific motivation, finally, which we repeat from the introduction, is that for some operators standardly used as conditionals in 3-valued logic, like the Strong Kleene conditional, no pure consequence relation (i.e. $ss$ or $tt$) can be such that this conditional satisfies the full deduction theorem relative to it. But a mixed relation such as $st$ does. Mixed consequence relations, which are relatively nonstandard, thus form a natural class to anyone concerned with the proof-theoretic behavior of operators that are standardly admitted.

\subsection{Order-theoretic consequence relations}

Before moving on, we also show that this approach to consequence relations subsumes another natural way to define a consequence relation, namely through an ordering on truth values:

\begin{definition}[Order theoretic relations]\label{def:ordertheoretic}
Suppose the set of truth values is equipped with an order $\leq$ {(a reflexive, transitive, antisymmetric relation{, be it partial or total})}, with $1$ ranked higher than any other value, and $0$ ranked lower. One may then define the \emph{order-theoretic consequence relation (associated to this order)}, requesting that premises have globally lower values than conclusions, through the formula:\medskip

\centerline{$\gamma\truthrelation_{\leq}\delta$ $\quad$ iff $\quad$ ($\exists x\in\gamma,\exists y\in\delta: x\leq y$) or ($0\in\gamma$ or $1\in\delta$).%
\footnote{One may entertain other ways to extend an order on truth values (possibly with more properties, such as the systematic presence of infimums and/or supremums) onto a truth-relation between subsets of truth values. Below are some examples, close to descriptions in \cite{chemla2017charac}, but which we will not attend to specifically here:
$\gamma\truthrelation\delta$ iff $\inf(\gamma)\leq\sup(\delta)$,
or 
$\gamma\truthrelation\delta$ iff $\exists d\in\delta: \inf(\gamma)\leq d$.
}}
\end{definition}

\noindent {We call truth-values other than 1 and 0 \emph{indeterminates}}. For illustration purposes, we distinguish two specific cases:
\begin{definition}[Total and degenerate order-theoretic relations]\label{def:specificordertheoretic}
$\phantom{x}$
\begin{itemize}
\item An order-theoretic relation is called \emph{total} if it is based on a total order of truth values: for all pairs of {indeterminates}, $\#_i\leq\#_j$ or $\#_j\leq\#_i$.
\item An order-theoretic relation is called \emph{degenerate} if it is based on the order such that no two {indeterminates} are {comparable}.

\end{itemize}
\end{definition}

We can here extend a result in the single-conclusion setting from \cite{chemla2017charac} to the current multi-conclusion framework:

\begin{theorem}\label{th:ordertheoric=intersectionpure}
An order-theoretic truth-relation derived from the order $\leq$ is equivalent to the intersection of all pure consequence relations based on designated values (containing $1$ but not $0$) that are upsets for this order.\footnote{Given an ordering $\leq$, an upset is a set that is closed under $\leq$, namely such that $y$ belongs to the set whenever $x$ belongs and $x\leq y$.}
\end{theorem}

\begin{myproof}
Choose $\gamma, \delta$. We must prove that the following statements are equivalent:
(i)~$\exists x\in\gamma, \exists y\in\delta: x\leq y$ or $0\in\gamma$ or $1\in\delta$ (that is, $\gamma\truthrelation_\leq\delta$),
and
(ii)~$\forall \mathcal{D} \textrm{ an upset of designated values}: \gamma\subseteq\mathcal{D}\Rightarrow\delta\cap\mathcal{D}\not=\emptyset$.

\begin{description}
\item[(i) entails (ii).] 
If $0\in\gamma$ or $1\in\delta$, the result follows because $\mathcal{D}s$ are sets of designated values (contain $1$ and not $0$). Otherwise, choose $(x_0, y_0)$ according to (i), and let $\mathcal{D}$ be an upset such that $\gamma\subseteq\mathcal{D}$. Then $x_0\in\gamma\subseteq\mathcal{D}$, so $y_0\in\mathcal{D}$, and therefore $y_0\in\delta\cap\mathcal{D}$, which is not empty.
\item[(ii) entails (i).] 
Let $\mathcal{D}=\{y | \exists x \in \gamma: x\leq y \}$. Then $\mathcal{D}$ is an upset (by transitivity of $\leq$) and $\gamma\subseteq\mathcal{D}$ (by reflexivity of $\leq$). Hence, if $\mathcal{D}$ is a set of designated values, then applying (ii) yields $y_0$ in $\delta\cap\mathcal{D}$, and it follows that there is $x_0\in\gamma$ such that $x_0\leq y_0$. If $\mathcal{D}$ is not a set of designated values, that would be either because $0\in\mathcal{D}$ or $1\not\in\mathcal{D}$. In the former case, that means that $0\in\gamma$. In the latter, it follows that $\gamma=\emptyset$, and so that $\gamma\subseteq\{1\}$, which is an upset of designated values and so, by (ii), that $\{1\}\cap\delta\not=\emptyset$, i.e.~$1\in\delta$.
\qedhere
\end{description}

\end{myproof}

This Theorem can be applied to specific examples:

\begin{corollary}\label{ex:ordertheoretictotal}
A total order-theoretic relation, with $0\leq\#_1\leq...\leq\#_N\leq 1$, is represented by:
$$\truthrelation_{\{1\},\{1\}}\cap\truthrelation_{\{1,\#_1\},\{1\,\#_1\}}\cap\truthrelation_{\{1,\#_1,\#_2\},\{1\,\#_1,\#_2\}}\cap...$$
\end{corollary}

\begin{corollary}\label{ex:ordertheoreticdegenerate}
The degenerate order-theoretic relation is represented by:
\begin{gather*}
	\truthrelation_{\{1\},\{1\}}	\\
	\cap \\
	\truthrelation_{\{1,\#_1\},\{1,\#_1\}} \cap \truthrelation_{\{1,\#_2\},\{1,\#_2\}} \cap ... \cap \truthrelation_{\{1,\#_N\},\{1,\#_N\}} 	\\
	\cap \\
	\truthrelation_{\{1,\#_1,\#_2\},\{1,\#_1,\#_2\}} \cap \truthrelation_{\{1,\#_1,\#_3\},\{1,\#_1,\#_3\}} ...  \cap \truthrelation_{\{1,\#_{N-1},\#_N\},\{1,\#_{N-1},\#_N\}}\\
	\cap\\ ...\\
	\cap \\
	\truthrelation_{\{1,\#_1,\#_2,..., \#_N\},\{1,\#_1,\#_2,..., \#_N\}}
\end{gather*}
\end{corollary}

The total and degenerate order-theoretic relations (and all other order-theoretic relations) in fact collapse in 3-valued logic, because there is no need to worry about the ordering among indeterminates when there is only one. The single order-theoretic relation in 3-valued logic is in fact one that has already been presented in Example~\ref{ex:sstt}, under the label $ss \cap tt$. 
With more than 3 truth-values, the notions come apart: an example of a degenerate (and non-total) order-theoretic relation is given by Belnap's 4-valued logic (\citealt{belnap1977useful}; see \citealt{shramko2011truth} for details on the order-theoretic aspect of the consequence relation, and \citealt{omori2015generalizing} for details on 4-valued connectives).
Order-theoretic relations form an important subclass of the intersective mixed consequence relations, and we will be able to provide extensive results concerning the types of connectives that they allow (see Theorem~\ref{th:fullordertheoretic}).

\subsection{Summary and assumptions}

We have here defined and delimited a space of consequence relations that we will explore systematically: unless otherwise noted, we will be interested in languages equipped with an intersective mixed consequence truth-relation and maximal or constant expressiveness.
In the next section, we turn to an examination of a particular type of constraints that can be put on logical connectives.

\section{Regular connectives\label{sec:operators}}

We are ultimately interested in specific connectives, such as conditionals, and their relation to consequence relations. Toward that goal, however, we shall take a broader view of other logical connectives, and ask what defines a connective from the point of view of a consequence relation. We want to find out how a given connective ought to interact with a consequence relation: what does it mean to have a formula headed by this connective in premise position? in conclusion position? These questions are of central importance in proof theory, and they were posed for the first time by \cite{gentzen1935investigations} when Gentzen introduced the sequent calculus for classical logic and for intuitionistic logic. We first review the way in which Gentzen's sequent rules arise for conjunction, disjunction, negation, and conditional, and then go on to introduce the class of what we call regular connectives for a consequence relation. We will provide various results about these connectives, revealing their behavior in classical logic, as well as their non-classical behavior in general. 
{We finally show that every Gentzen connective can be exemplified by a truth-function in a particular $4$-valued logic.}

\subsection{Gentzen's operational rules, first examples}\label{sec:gentzen}

Gentzen's sequent calculus for classical logic can be viewed as a systematic framework to represent the way in which the meaning of a logical connective ought to be understood in relation to a given consequence relation. Gentzen distinguished \emph{structural rules}, concerned with the general behavior of the consequence relation (is it contractive, commutative, monotone, etc), from \emph{operational rules}, concerned with the specific behavior of logical connectives. For each connective, Gentzen's operational rules specify how the connective is to be treated as a premise in an argument, or as a conclusion in an argument. 

Consider conjunction and disjunction first. In Gentzen's approach, the concatenation of premises corresponds to their conjunction, and the concatenation of conclusions correspond to their disjunction. This yields the following rules when a conjunctive formula appears as a premise, and a disjunctive formula as a conclusion:

\begin{itemize}
\item $\forall\Gamma,\Delta: \Gamma, P \wedge Q \vdash \Delta$ iff $\Gamma, P, Q \vdash \Delta$
\item $\forall\Gamma,\Delta: \Gamma \vdash P \vee Q, \Delta$ iff $\Gamma, \vdash P, Q, \Delta$
\end{itemize}

\noindent When a conjunctive formula appears as a conclusion, or a disjunctive formula as a premise, we get the following dual rules:
\begin{itemize}
\item $\forall\Gamma,\Delta: \Gamma \vdash P \wedge Q, \Delta$ iff $\Gamma \vdash P, \Delta$ and $\Gamma \vdash Q, \Delta$
\item $\forall\Gamma,\Delta: \Gamma, P \vee Q \vdash \Delta$ iff $\Gamma, P \vdash \Delta$ and $\Gamma, Q \vdash \Delta$
\end{itemize}

For negation, we get the following sequent rules for when negation appears in premise position, or in conclusion position:

\begin{itemize}
\item $\forall\Gamma,\Delta: \Gamma, \neg P \vdash \Delta$ iff $\Gamma \vdash P, \Delta$
\item $\forall\Gamma,\Delta: \Gamma \vdash \neg P, \Delta$ iff $\Gamma, P \vdash \Delta$
\end{itemize}

For the conditional {(or implication)}, Gentzen proposed the following rules:\footnote{Gentzen originally stated only the right-to-left direction of those rules, but it is natural to use invertible rules.}
\begin{itemize}
\item $\forall\Gamma,\Delta: \Gamma \vdash P \rightarrow Q, \Delta$ 
	iff $\Gamma, P \vdash Q, \Delta$
\item $\forall\Gamma,\Delta: \Gamma , P \rightarrow Q\vdash \Delta$ 
	iff ($\Gamma \vdash P, \Delta$ and $\Gamma, Q \vdash \Delta$)
\end{itemize}

The first of these rules basically corresponds to the full deduction theorem: it tells us that the conditional internalizes the consequence relation in the object-language. The rule for the conditional in premise position is best understood in relation to the standard definition of logical consequence in classical logic: a conditional is not designated (false) provided its antecedent is designated and its consequent is not.

\subsection{Regularity: definition and first application}

An important feature of the Gentzen rules is their analytic character, namely the fact that the meaning of a connective, whether in premise position or in conclusion position, is explained fully in terms of the (possibly empty) conjunction of sequent rules involving the component subformulae of the formula build from that connective. In \cite{chemla2018suszko}, we put forward the notion of a \emph{regular connective} to describe logical connectives obeying such constraints:

\begin{definition}[Gentzen-regular connectives]\label{def:regconn}
	Given a consequence relation $\vdash$, an $n$-ary connective $C$ is \emph{regular} for it if there exist
	$\mathcal{B}^p\subseteq\mathcal{P}(\{1,..., n\})\times\mathcal{P}(\{1,..., n\})$ and
	$\mathcal{B}^c\subseteq\mathcal{P}(\{1,..., n\})\times\mathcal{P}(\{1,..., n\})$
such that
	$\forall\Gamma, \Delta, \forall F_1, ..., F_n:$
	\[\begin{array}{c@{\textrm{ iff }}c}
	\Gamma, C(F_1, ..., F_n) \vdash \Delta
		& \bigwedge\limits_{(B_p,B_c)\in \mathcal{B}^p} 
			{\Gamma, \{F_i: i\in B_p\}\vdash \{F_i: i\in B_c\}, \Delta}\\

	\Gamma \vdash C(F_1, ..., F_n) , \Delta
		& \bigwedge\limits_{(B_p,B_c)\in \mathcal{B}^c} 
			{\Gamma, \{F_i: i\in B_p\}\vdash \{F_i: i\in B_c\}, \Delta}\\
	\end{array}\]

		\end{definition}

To immediately provide examples, we state that every truth-function in classical logic defines a regular connective. We simply provide the regularity rules associated with binary truth-functions, we will prove the general result later in Section~\ref{sec:bivalentregconnectives}.

\begin{example}[Regular rules associated to all binary truth-functional connectives in classical logic]\label{regrulesclassicalconnectives} All 16 binary truth-functions of classical logic can be associated with a regularity rule, as indicated in the following table.\footnote{This is not to imply that a regular connective may satisfy only one regularity rule. For instance, in reflexive logics, as in classical logic, adding a conjunct of the form $\Gamma, A \vdash A, \Delta$ to a regularity rule produces a new rule, but it is essentially the same rule and certainly it is satisfied by the same connectives.}
\medskip
	
\centering\noindent {\small\framebox{\begin{tabular}{r|l|ll}
	
$O(P,Q)$		& Premise-regularity rules & Conclusion-regularity rules\\
		& $\Gamma, O(P,Q) \vdash \Delta$ iff: & $\Gamma \vdash O(P,Q), \Delta$ iff:\\
	
	\hline
	\hline
	
	$\top$
		& $\Gamma \vdash \Delta$
		& \emph{always true} (an empty conjunction)\\
	$\bot$
		& \emph{always true} (an empty conjunction)
		& $\Gamma \vdash \Delta$\\

	\hline
	
	$P$
		& $\Gamma, P \vdash \Delta$
		& $\Gamma \vdash P, \Delta$\\

	$\neg P$
		& $\Gamma \vdash P, \Delta$
		& $\Gamma, P \vdash \Delta$\\

	\hline

	$Q$
		& $\Gamma, Q \vdash \Delta$
		& $\Gamma \vdash Q, \Delta$\\

	$\neg Q$
		& $\Gamma \vdash Q, \Delta$
		& $\Gamma, Q \vdash \Delta$\\

	\hline

	$P\vee Q$
		& $\Gamma, P \vdash \Delta$ and $\Gamma, Q \vdash \Delta$
		& $\Gamma \vdash P, Q, \Delta$\\

	$\neg (P\vee Q)$
		& $\Gamma \vdash P, Q, \Delta$
		& $\Gamma, P \vdash \Delta$ and $\Gamma, Q \vdash \Delta$\\

	\hline

	$P\wedge Q$
		& $\Gamma, P, Q \vdash \Delta$
		& $\Gamma \vdash P, \Delta$ and $\Gamma \vdash Q, \Delta$\\

	$\neg (P\wedge Q)$ 
		& $\Gamma \vdash P, \Delta$ and $\Gamma \vdash Q, \Delta$
		& $\Gamma, P, Q \vdash \Delta$\\

	\hline

	$P\rightarrow Q$
		&  $\Gamma, Q \vdash \Delta$ and $\Gamma \vdash P, \Delta$ 
		& $\Gamma, P \vdash Q, \Delta$\\

	$P\wedge \neg Q$ 
		& $\Gamma, P \vdash Q, \Delta$
		& $\Gamma, Q \vdash \Delta$ and $\Gamma \vdash P, \Delta$\\

	\hline

	$P\leftarrow Q$
		& $\Gamma, P \vdash \Delta$ and $\Gamma \vdash Q, \Delta$ 
		& $\Gamma, Q \vdash P, \Delta$\\

	$\neg P\wedge Q$ 
		& $\Gamma, Q \vdash P, \Delta$
		& $\Gamma, P \vdash \Delta$ and $\Gamma \vdash Q, \Delta$\\

	\hline

	$P\leftrightarrow Q$ 
		& $\Gamma, P, Q \vdash \Delta$ and $\Gamma \vdash P, Q, \Delta$
		& $\Gamma, P \vdash Q, \Delta$ and $\Gamma, Q \vdash P, \Delta$\\

	$P~\underline{\vee}~Q$  
		& $\Gamma, P \vdash Q, \Delta$ and $\Gamma, Q \vdash P, \Delta$
		& $\Gamma, P, Q \vdash \Delta$ and $\Gamma \vdash P, Q, \Delta$ 
		\\

\end{tabular}}}\bigskip

\end{example}

\subsection{Interactions between G-connectives: some classical rules}

We state here a few properties concerning the interactions between Gentzen-regular connectives. These highlight, again, the similarity of behavior with classical, truth-functional connectives. We start with a general closure condition.

\begin{theorem}
A combination of Gentzen-regular connectives is Gentzen-regular.
\end{theorem}

\begin{myproof}
This can be proved by induction, a conjunction of conjunctions being a conjunction itself.
\end{myproof}

\noindent 
Furthermore, every $n$-ary projection $pr_i^n$ is Gentzen regular, which together with the previous result establishes that Gentzen regular connectives form what is called a \emph{clone} (see \citealp{KERKHOFF2014107} for a full presentation).

\begin{theorem}
A projection $pr_i^n$, which associates $pr_i^n(X_1,..., X_n)$ to $X_i$, satisfies the regularity rules:
\begin{itemize}
\item $\Gamma, pr_i^n(X_1,..., X_n) \vdash \Delta \textrm{ iff } \Gamma, X_i \vdash \Delta$
\item $\Gamma \vdash pr_i^n(X_1,..., X_n), \Delta \textrm{ iff } \Gamma \vdash X_i, \Delta$
\end{itemize}
\end{theorem}

\begin{myproof}
Immediate.
\end{myproof}

Let us define now a couple of Gentzen-regular connectives of particular salience:

\begin{definition}[Gentzen-connectives]

We call a G-conjunction / G-disjunction /  G-negation / G-conditional a truth-functional connective which obeys the appropriate Gentzen rules stated in Section~\ref{sec:gentzen}.

\end{definition}

\noindent The interactions between these operators are classical. For instance, the usual De Morgan's laws and the contraposition rule hold:

\begin{theorem}\label{th:demorgan}[De Morgan's laws]
$\phantom{x}$
\begin{itemize}
\item  If $\neg$ is a G-negation and $\wedge$ is a G-conjunction, then $\neg(\neg\_\_\wedge\neg\_\_)$ is a G-disjunction.
\item  If $\neg$ is a G-negation and $\vee$ is a G-disjunction, then $\neg(\neg\_\_\vee\neg\_\_)$ is a G-conjunction.
\end{itemize}
\end{theorem}

\begin{theorem}\label{th:contraposition}(Contraposition)
If $\to$ is a G-conditional and $\neg$ a G-negation, then the contraposition rules apply: 
\begin{itemize}
\item $\Gamma, A\to B \vdash \Delta$ iff $\Gamma, (\neg B \to \neg A) \vdash \Delta$
\item $\Gamma \vdash A\to B, \Delta$ iff $\Gamma \vdash (\neg B \to \neg A),  \Delta$
\end{itemize}
\end{theorem}

\begin{myproof}[Theorem~\ref{th:demorgan}]
For instance, we can show the premise-regularity rule for the G-disjunction defined from a G-negation and a G-conjunction:
$\Gamma, \neg(\neg A \wedge\neg B)\vdash\Delta$
iff
	$\Gamma \vdash (\neg A \wedge\neg B), \Delta$ 
	(negation premise-regularity rule)
iff
	$\Gamma \vdash \Delta$ and $\Gamma \vdash \neg B, \Delta$ 
	(conjunction conclusion-regularity rule)
iff
	$\Gamma, A \vdash \Delta$ and $\Gamma, B \vdash \Delta$ 
	(negation conclusion-regularity rule twice), which is the disjunction premise-regularity rule.
\end{myproof}

\begin{myproof}[Theorem~\ref{th:contraposition}]
\begin{itemize}
\item $\Gamma, (\neg B \to \neg A) \vdash \Delta$ 
iff $\Gamma,  \neg A \vdash \Delta$ and $\Gamma \vdash \neg B, \Delta$  (premise Gentzen-regularity rule for $\to$)
iff $\Gamma \vdash A, \Delta$ and $\Gamma, B \vdash \Delta$  (Gentzen-regularity rules for $\neg$) 
iff $\Gamma,  (A \to B) \vdash \Delta$  (premise Gentzen-regularity rule for $\to$).
\item $\Gamma \vdash (\neg B \to \neg A) , \Delta$ 
iff $\Gamma,  \neg B \vdash \neg A, \Delta$  (conclusion Gentzen-regularity rule for $\to$)
iff $\Gamma,  A \vdash B, \Delta$  (Gentzen-regularity rules for $\neg$) 
iff $\Gamma \vdash (A \to B), \Delta$  (conclusion Gentzen-regularity rule for $\to$).
\qedhere\end{itemize}
\end{myproof}

Given that G-conditionals are our core interest, we here introduce a couple of specific results. We show in particular that the existence of a G-conditional may be guaranteed by the existence of other G-connectives, and that it actually guarantees the existence of the other G-connectives assuming that there is a formula $\bot$ with constant value $0$ (which is one of our default assumptions for the language and its semantics, see Definition~\ref{def:constantsandmaximalexpressiveness}). Theorem~\ref{th:conditionalissufficient} will later provide a more general result of complete expressiveness based on G-conditionals however, within section~\ref{sec:bivalentregconnectives} which discusses bivalent connectives more systematically.

\begin{theorem}[Interactions between G-conditionals and other G-connectives]\label{th:combinationsforconditionals}
$\phantom{x}$
\begin{itemize}
\item  If $\neg$ is a G-negation and $\vee$ is a G-disjunction, then $(\neg\_\_)\vee\_\_$ is a G-conditional,
\item  If $\neg$ is a G-negation and $\wedge$ is a G-conjunction, then $\neg(\_\_\wedge\neg\_\_)$ is a G-conditional,
\item  If $\to$ is a G-conditional, $\_\_\to\bot$ is a G-negation,
\item  If $\to$ is a G-conditional, $(\_\_\to\bot)\to \_\_$ is a G-disjunction,
\item  If $\to$ is a G-conditional,  $(\_\_\to (\_\_\to\bot))\to\bot$ is a G-conjunction.
\end{itemize}
\end{theorem}

\begin{myproof}
The proof is essentially the same as before and obtained by combining regularity rules with one another, here using in particular the regularity rule of $\bot$ as given in Example~\ref{regrulesclassicalconnectives}.
\end{myproof}

\subsection{Non-classical behavior}

Crucially, G-regularity does not entail a completely classic behavior, unless the consequence relation itself plays its part: for instance, a number of classical rules will hold with G-regular connectives if and only if the consequence relation is reflexive.

\begin{theorem}[Reflexivity and G-connectives]
Consider a language with a monotonic
consequence relation, a constant expressive semantics, and G-connectives $\neg$, $\vee$, $\wedge$ and $\to$. The following statements are all equivalent:
\begin{enumerate}
\item $\forall A:\ A \vdash A$
	\hfill \emph{(Reflexivity)}
\item $\forall A, B:\ (A \to B) \vdash (A \to B)$ 
	\hfill \emph{(Reflexivity restricted to conditional formulae)}
\item $\forall A:\ \neg A \vdash \neg A$ 
	\hfill \emph{(Reflexivity restricted to negative formulae)}
\item $\forall A, B:\ (A \wedge B) \vdash (A \wedge B)$
	\hfill \emph{(Reflexivity restricted to conjunctive formulae)}
\item $\forall A, B:\ (A \vee B) \vdash (A \vee B)$
	\hfill \emph{(Reflexivity restricted to disjunctive formulae)}

\item $\forall A, B:\ A, (A \to B) \vdash B$ 
	\hfill \emph{(Modus Ponens)}

\item $\forall A:\ \vdash A, \neg A$
	\hfill \emph{(Law of Excluded Middle, version 1)}
\item $\forall A:\ \vdash (A \vee \neg A)$
	\hfill \emph{(Law of Excluded Middle, version 2)}

\item $\forall A:\ A, \neg A \vdash $
	\hfill \emph{(Principle of Explosion, version 1)}
\item $\forall A:\ (A \wedge \neg A)\vdash$
	\hfill \emph{(Principle of Explosion, version 2)}

\item $\forall A, B:\ A \vdash A, B$ 
	\hfill \emph{(Intermediate step)}
\item $\forall A, B:\ A, B \vdash A$ 
	\hfill \emph{(Intermediate step)}
\item $\forall A, B:\ A, B \vdash A, B$ 
	\hfill \emph{(Intermediate step)}

\end{enumerate}

\end{theorem}

\begin{myproof}
We start by showing the equivalence between 1, 11, 12 and 13: By monotonicity, 1 entails 11, 12 and 13; conversely, if we choose $B=A$ in either 11, 12 or 13, we see that either of them entails 1. From there, all the rest follows, because every statement from 2 to 10 all can be reduced to a conjunction of statements among 1, 11, 12, 13 by applying the regularity rule of the relevant connective whenever it appears. 

Let us take three examples. First, Modus Ponens (6):
$\forall A, B:\ A, (A \to B) \vdash B$ 
is true iff
$\forall A, B:\ A, B \vdash A, B$, by applying the premise regularity rule of the G-conditional.

Second, the Law of Excluded Middle, version 2 (8):
$\forall A:\ \vdash (A \vee \neg A)$. There is a G-disjunction in the conclusions, so we apply the conclusion regularity for a G-disjunction and obtain that this is true iff
$\forall A:\ \vdash A, \neg A$. This is (7). Since it contains a G-negation in conclusion, we can apply the relevant rule and obtain that this is equivalent to $\forall A:\ A \vdash A$.

Finally, we can show why reflexivity restricted to conditional formulae (2) is equivalent to plain reflexivity:
$\forall A, B:\ (A \to B) \vdash (A \to B)$ 
iff
$\forall A, B:\ (A \to B), A \vdash B$  (conclusion regularity rule for the G-conditional)
iff
$\forall A, B:\ B, A \vdash B \textrm{ and } A \vdash A, B$  (premise regularity rule for the G-conditional).
iff
$\forall A, B:\ B, A \vdash B \textrm{ and } \forall A, B: A \vdash A, B$  (distributivity of conjunction and universal quantifiers).
And this last statement is indeed a conjunction of the statements 12 and 13 (which are, remember, equivalent).
\end{myproof}

Not all intersective mixed consequence relations are reflexive, so this result reveals that Gentzen-regular rules are not sufficient to impose a classical behavior for the connectives if the consequence relation itself is not structurally classical in the first place. Similarly, we know that not all intersective mixed consequence relations are transitive. Hence, we could investigate what laws of classical logic are equivalent to transitivity, or any other such property for that matter. Instead we now turn to the issue of whether G-connectives always exist for a consequence relation.

\subsection{Regularity: key results}

This section presents more consequences of regularity. These properties are fundamental, and they will help to prove important results later on. Some of them are rather technical and can be skipped on a first reading.

To begin with, in order to investigate the connection between regular connectives and truth-relational consequence relations, it is useful to assume constant expressiveness. Under that assumption, regular connectives can be translated at the level of truth values without any loss.

\begin{theorem}\label{th:reducetotruthvalues}
In a language with maximal or constant expressiveness (see Definition~\ref{def:constantsandmaximalexpressiveness}), a truth-functional connective $C$ with truth-function $\underline{C}$ satisfies a regularity rule iff it satisfies it at the level of truth values:\medskip

\begin{tabular}{lll}
	&
	$\forall\Gamma, \Delta, \forall{F_1, ..., F_n}:$ &
	$\begin{array}[t]{c@{\quad\textrm{ iff }\quad}c}
	\Gamma, C(F_1, ..., F_n) \vdash \Delta
		& \bigwedge\limits_{(B_p,B_c)\in \mathcal{B}^p} 
			{\Gamma, \{F_i: i\in B_p\}\vdash \{F_i: i\in B_c\}, \Delta}\\

	\Gamma \vdash C(F_1, ..., F_n) , \Delta
		& \bigwedge\limits_{(B_p,B_c)\in \mathcal{B}^c} 
			{\Gamma, \{F_i: i\in B_p\}\vdash \{F_i: i\in B_c\}, \Delta}\\
	\end{array}$
	
	\\

iff \\

	&
	$\forall\gamma, \delta, \forall x_1, ..., x_n:$ &
	$\begin{array}[t]{c@{\quad\textrm{ iff }\quad}c}
	\gamma, \underline{C}(x_1, ..., x_n) \truthrelation \delta
		& \bigwedge\limits_{(B_p,B_c)\in \mathcal{B}^p} 
			{\gamma, \{x_i: i\in B_p\}\truthrelation \{x_i: i\in B_c\}, \delta}\\

	\gamma \truthrelation \underline{C}(x_1, ..., x_n) , \delta
		& \bigwedge\limits_{(B_p,B_c)\in \mathcal{B}^c} 
			{\gamma, \{x_i: i\in B_p\}\truthrelation \{x_i: i\in B_c\}, \delta}\\
	\end{array}$
	\end{tabular}

\end{theorem}

\begin{myproof}
The downward to upward direction is clear: if the rule holds for every possible list of truth values, it holds for every semantic interpretation taken as a whole. Conversely, because the language is constant expressive, the fact that the rule holds for \emph{all} formulae, therefore for all constant formulae, guarantees that it holds in the second form too (when the formulae are constant, the universal statement that reduces $\vdash$ to $\truthrelation$ is trivial).
\end{myproof}

Secondly, regularity and the truth-relation together put very explicit constraints on the truth-function of a regular connective, of which we can extract a useful constructive presentation:

\begin{theorem}\label{th:truthconstraintforregconn}
Consider a logic equipped with a truth-relation with a minimal representation
$\truthrelation{\mathcal{D}_p^1,\mathcal{D}_c^1}\cap...\cap\truthrelation{\mathcal{D}_p^K,\mathcal{D}_c^K}$, and an $n$-ary connective $C$ with truth-function $\underline{C}$. $C$ follows a regularity rule (defined as above) iff
$\forall k, \forall x_1, ..., x_n$:

	\[\begin{array}{c@{\quad\textrm{ iff }\quad}c}

	\underline{C}(x_1, ..., x_n) \not\in\mathcal{D}_p^k
		& 
		\bigwedge\limits_{(B_p,B_c)\in \mathcal{B}^p} 
			{\{x_i: i\in B_p\}\subseteq\mathcal{D}_p^k \Rightarrow \{x_i: i\in B_c\}\cap\mathcal{D}_c^k\not=\emptyset}
			\\

	\underline{C}(x_1, ..., x_n) \in\mathcal{D}_c^k
		& 
		\bigwedge\limits_{(B_p,B_c)\in \mathcal{B}^c} 
			{\{x_i: i\in B_p\}\subseteq\mathcal{D}_p^k \Rightarrow  \{x_i: i\in B_c\}\cap\mathcal{D}_c^k\not=\emptyset}
			\\

	\end{array}\]
\end{theorem}

\begin{myproof}
Consider a connective satisfying these constraints, and choose $\gamma, \delta, x_1, ..., x_n$.

$\gamma, \underline{C}(x_1,...,x_n)\truthrelation\delta$
\begin{minipage}[t]{.8\textwidth}
iff
$\forall k:\gamma, \underline{C}(x_1,...,x_n)\truthrelation_k\delta$

iff
$\forall k: [\underline{C}(x_1,...,x_n)\not\in\mathcal{D}_p^k
	\textrm{ or }
	\gamma \truthrelation_k \delta]$

iff
$\forall k: [
		(\bigwedge\limits_{(B_p,B_c)\in \mathcal{B}^p} 
			{\{x_i: i\in B_p\}\subseteq\mathcal{D}_p^k \Rightarrow \{x_i: i\in B_c\}\cap\mathcal{D}_c^k\not=\emptyset})
	\textrm{ or }
	\gamma \truthrelation_k \delta]$

iff
$\forall k: 
		\bigwedge\limits_{(B_p,B_c)\in \mathcal{B}^p} 
			{\{x_i: i\in B_p\}, \gamma \truthrelation_k \{x_i: i\in B_c\}, \delta}$

iff
$\bigwedge\limits_{(B_p,B_c)\in \mathcal{B}^p} 
			{\{x_i: i\in B_p\}, \gamma \truthrelation \{x_i: i\in B_c\}, \delta}$.
\end{minipage}

\noindent This proves that the premise regularity rule holds, and an analogous derivation would prove that the conclusion regularity rule also holds.

Conversely, consider a connective satisfying the relevant regularity rule. Then pick $k_0, x_1, ..., x_n$. First note that 
$\forall k\not=k_0: 
\mathcal{D}_p^{k_0}\truthrelation_{\mathcal{D}_p^k,\mathcal{D}_c^k} \overline{\mathcal{D}_c^{k_0}}$, for otherwise 
$\mathcal{D}_p^{k_0}\subseteq\mathcal{D}_p^k$ and 
$\mathcal{D}_c^k\subseteq\mathcal{D}_c^{k_0}$, in which case $\truthrelation_{k_0}$ should have been dropped from the minimal representation of $\truthrelation$. 
Then consider the statement $(\mathcal{S}): \mathcal{D}_p^{k_0}, C(x_1, ..., x_n)\truthrelation \overline{\mathcal{D}_c^{k_0}}$.
\begin{itemize}
\item On the one hand, given that $\forall k\not=k_0: \mathcal{D}_p^{k_0}\truthrelation_k \overline{\mathcal{D}_c^{k_0}}$ is already true,
$(\mathcal{S})$ is true iff
$\mathcal{D}_p^{k_0}, C(x_1, ..., x_n)\truthrelation_{k_0} \overline{\mathcal{D}_c^{k_0}}$, which boils down to whether $C(x_1, ..., x_n)\not\in\mathcal{D}_p^{k_0}$.
\item On the other hand, by regularity, $(\mathcal{S})$ is true 
	\begin{minipage}[t]{.6\textwidth}
		iff
		$\bigwedge\limits_{(B_p,B_c)\in \mathcal{B}^p} 
			{\mathcal{D}_p^{k_0}, \{x_i: i\in B_p\}\truthrelation \{x_i: i\in B_c\}, \overline{\mathcal{D}_c^{k_0}}}$
			
		iff 
		$\forall k: \bigwedge\limits_{(B_p,B_c)\in \mathcal{B}^p} 
			{\mathcal{D}_p^{k_0}, \{x_i: i\in B_p\}\truthrelation_k \{x_i: i\in B_c\}, \overline{\mathcal{D}_c^{k_0}}}$

		iff 
		$\bigwedge\limits_{(B_p,B_c)\in \mathcal{B}^p} 
			{\mathcal{D}_p^{k_0}, \{x_i: i\in B_p\}\truthrelation_{k_0} \{x_i: i\in B_c\}, \overline{\mathcal{D}_c^{k_0}}}$

		iff 
			$\{x_i: i\in B_p\}\subseteq\mathcal{D}_p^{k_0} \Rightarrow \{x_i: i\in B_c\}\cap\mathcal{D}_c^{k_0}\not=\emptyset$.
	\end{minipage}
\end{itemize}
\noindent
Putting together the final step of the two equivalence derivations above provides the first half of the result, corresponding to the premise regularity rule (first of the two conditions). The second half is obtained similarly by manipulating the statement 
$(\mathcal{S'}): \mathcal{D}_p^{k_0} \truthrelation C(x_1, ..., x_n), \overline{\mathcal{D}_c^{k_0}}$.
\end{myproof}

Corollary~\ref{th:formcond} illustrates how this Theorem~\ref{th:truthconstraintforregconn} constrains the form of G-conjunctions, G-disjunctions and G-conditionals. A useful paraphrase of the regularity rules for the conditional emerges: when the value of the first argument of a conditional is designated, then the value of the second argument must be designated too (abstracting away from the difference between premise and conclusion designated values, or from the various mixed consequence relations potentially involved). This paraphrase surely is reminiscent of the intuition behind the definition of a mixed consequence relation in the first place.

\begin{corollary}\label{th:formcond}
Consider a truth-relation $\truthrelation$ minimally represented as $\truthrelation_{\mathcal{D}_p^1, \mathcal{D}_c^1}\cap...\cap\truthrelation_{\mathcal{D}_p^K, \mathcal{D}_c^K}$.
\begin{itemize}
\item A connective $\wedge$ is a G-conjunction for $\truthrelation$ iff $\forall a, b$, for all sets of designated values $\mathcal{D}$: %(premise or conclusion):
\begin{itemize}
\item $(a \wedge b) \in \mathcal{D} \quad$ iff $\quad$ ($a\in\mathcal{D}  \textrm{ and }  b\in\mathcal{D}$).
\end{itemize}

\item A connective $\vee$ is a G-disjunction for $\truthrelation$ iff $\forall a, b$, for all sets of designated values $\mathcal{D}$: %(premise or conclusion):
\begin{itemize}
\item $(a \vee b) \in \mathcal{D} \quad$ iff $\quad$ ($a\in\mathcal{D} \textrm{ or } b\in\mathcal{D}$).
\end{itemize}

\item A connective $\to$ is a G-conditional for $\truthrelation$ iff $\forall a, b:\forall i:$
\begin{itemize}
\item $(a \to b) \in \mathcal{D}_p^i \quad$ iff $\quad$ ($a\in\mathcal{D}_c^i \Rightarrow b\in\mathcal{D}_p^i$)
\item $(a \to b) \in \mathcal{D}_c^i \quad$ iff $\quad$ ($a\in\mathcal{D}_p^i \Rightarrow b\in\mathcal{D}_c^i$)
\end{itemize}
\end{itemize}
\end{corollary}

In section~\ref{sec:N}, we will capitalize on these constraints on G-connectives to characterize what truth-relations can admit a truth-function for them. For now, we can mention a more general consequence of this Theorem: the regular connectives of an \emph{intersective} mixed consequence truth-relation are inherited from the regular connectives of the mixed consequence truth-relations they are made of. Consider an intersective mixed consequence truth-relation, say $ss\cap tt$. A connective $C$ can be a regular connective for $ss \cap tt$, say a G-conjunction, if and only if $\underline{C}$ is a common G-conjunction of $ss$ and of $tt$ in the first place. More generally:

\begin{corollary}\label{cor:intersectioncoincideR}
Consider an intersective mixed consequence truth-relation with a minimal representation written as $\truthrelation_{\mathcal{D}_p^1,\mathcal{D}_c^1}\cap...\cap\truthrelation_{\mathcal{D}_p^K,\mathcal{D}_c^K}$. And consider some regularity rule $\mathcal{R}$. Then
\begin{tabular}{ll}
& $f$ is the truth-function of a regular connective satisfying $\mathcal{R}$ for $\truthrelation$ \\
iff &  $f$ is the truth-function of a regular connective satisfying $\mathcal{R}$ for each $\truthrelation_{\mathcal{D}_p^k,\mathcal{D}_c^k}$.
\end{tabular}
\end{corollary}

In a similar note, we can consider the reverse question: given some operator, say a conditional operator, what consequence relation may be associated to it? Concretely, the following result entails in a compact logic {(see Definition \ref{def:regularcompactlogic})} that $ss$ and $ts$ have no conditional in common (because they have the same `conclusion-tautologies' ($\{1\}$)); and that $ss$ and $st$ neither (because they have the same `premise-contradictions' ($\{0,\half\}$)).%\nbPE{does it really matter to mention compactness here? The proof does not refer to compactness explicitly below, worth clarifying in my opinion}

\begin{theorem}\label{th:fromcondtoLC}
Two truth-relational consequence relations $\vdash_1$ and $\vdash_2$ in a compact logic are identical as soon as they share (i)~a truth-functional conditional and (ii)~the same conclusion-tautologies ($\{x: \vdash_i x\}$) or the same premise-contradictions ($\{x: x \vdash_i\}$).
\end{theorem}

\begin{myproof}
Assume that $\to$ is a conditional for $\vdash_1$ and $\vdash_2$, then the two relations also share a G-conjunction, a G-disjunction and a G-negation. For a given finite argument $P_1, ..., P_n\vdash Q_1, ..., Q_n$, all premises can be moved one by one to the conclusions, using the regularity rule for the G-negation, hence the argument holds iff 
$\vdash \neg P_1, ..., \neg P_n, Q_1, ..., Q_n$. 
We can continue and collapse all conclusions with the regularity rule for disjunction, so the argument holds iff
$\vdash \neg P_1 \vee ... \vee \neg P_n \vee Q_1 \vee ... \vee Q_n$. 
Hence, given G-disjunctions and G-negations, whether {a finite} argument holds only depends on what {arguments} of the form {$\vdash Q$} {hold}. If two consequence relations share a G-disjunction and a G-negation, and have the same conclusion tautologies, they are identical {on finite arguments, and if the logic is compact they are thus identical}.

Similarly, if two consequence relations share a G-conjunction and G-negation, and have the same premise-contradictions, they are identical because every argument $P_1, ..., P_n\vdash Q_1, ..., Q_n$ boils down to
$P_1 \wedge ... \wedge P_n \wedge \neg Q_1 \wedge ... \wedge \neg Q_n \vdash $. 
\end{myproof}

To conclude this section, we explore the (limited) multiplicity of regular connectives of any given type, for a given truth-relation.
Whether a truth-relation admits a connective satisfying a particular regularity rule is in general a difficult question. Before diving into it in the following sections, we address here a different problem: \emph{given a consequence relation and a regularity rule, how many truth-functions can produce a connective satisfying this rule, if any?}
The solution is mostly guided by the truth-relation, not by the specifics of the regularity rule. 
Theorem \ref{th:systematicmulitplicity} states conditions under which a truth-relation allows for a multiplicity of regular connectives. The converse result (Theorem~\ref{th:limitedmulitplicity}) shows that this is the total extent to which there could be multiplicity.

\begin{theorem}\label{th:systematicmulitplicity}
Suppose two truth values $y_1$ and $y_2$ play the same role throughout the truth-relation (that is, 
$\forall\gamma, \delta: \gamma, y_1\truthrelation \delta$ iff $\gamma, y_2\truthrelation \delta$
and
$\forall\gamma, \delta: \gamma \truthrelation y_1, \delta$ iff $\gamma \truthrelation  y_2, \delta$). If $f$ is a truth-function for a connective that follows some regularity rule and such that $f(\vec{x})=y_1$, then the truth-function $f_{\vec{x}\to{y_2}}$ which is just like $f$ except that $f_{\vec{x}\to{y_2}}(\vec{x})=y_2$ follows the same regularity rule.
\end{theorem}

\begin{myproof}
Consider a connective $C$ to which the truth-function $f$ is assigned and which satisfies some regularity rule, for instance a premise-regularity rule:
$\Gamma, C(A_1,..., A_n)\vdash\Delta$ 
iff 
$\bigwedge_{(B_p,B_c)\in\mathcal{B}^p} \Gamma, \{F_i : i \in \mathcal{B}_p\} \vdash\{F_i : i \in \mathcal{B}_c\}, \Delta$.
Now we can show that a connective $C'$ which has $f_{\vec{x}\to{y_2}}$ has a truth-function will follow the same regularity rule. 
$\Gamma, C'(A_1,..., A_n)\vdash\Delta$ 
iff 
$\forall v: v(\Gamma), v(C')(v(A_1),..., v(A_n))\truthrelation v(\Delta)$  (by definition of a truth-relation)
iff 
$\forall v: v(\Gamma), f_{\vec{x}\to{y_2}}(v(A_1),..., v(A_n))\truthrelation v(\Delta)$  (applying truth-functionality)
iff 
$\forall v: v(\Gamma), f(v(A_1),..., v(A_n))\truthrelation v(\Delta)$ (because the difference between $f$ and $f_{\vec{x}\to{y_2}}$ is irrelevant to such arguments)
iff 
$\forall v: v(\Gamma), v(C)(v(A_1),..., v(A_n))\truthrelation v(\Delta)$ (by truth-functionality of $C$)
iff 
$\Gamma, C(A_1,..., A_n)\vdash\Delta$   (by definition of the truth-relation)
iff
$\bigwedge_{(B_p,B_c)\in\mathcal{B}^p} \Gamma, \{F_i : i \in \mathcal{B}_p\} \vdash\{F_i : i \in \mathcal{B}_c\}, \Delta$   (applying the regularity rule to $C$).
\end{myproof}

\begin{theorem}[Limited multiplicity]\label{th:limitedmulitplicity}
Consider $f_1$ and $f_2$, two truth-functions associated with a given regularity rule. Suppose that they differ on the value they assign to some input $\vec{x}$: $f_1(\vec{x})\not=f_2(\vec{x})$. Then these two values $f_1(\vec{x})$ and $f_2(\vec{x})$ play the same role throughout the truth-relation.
\end{theorem}

\begin{myproof}
The following two equivalences are obtained by plugging in each time the relevant regularity rule:
$\forall \gamma, \delta: \gamma, f_1(\vec{x}) \truthrelation \delta$ iff $\gamma, f_2(\vec{x}) \truthrelation \delta$
and
$\forall \gamma, \delta: \gamma \truthrelation f_1(\vec{x}), \delta$ iff $\gamma \truthrelation f_2(\vec{x}), \delta$.
\end{myproof}

\subsection{Regularity: classical logic, bivalent connectives}\label{sec:bivalentregconnectives}

With the previous technical results in place, we can prove a result announced earlier:

\begin{theorem}\label{th:classicalallregular}
Every truth-functional connective in classical logic is regular.
\end{theorem}

\begin{myproof}
Suppose $C$ is an $n$-ary truth-functional connective in classical logic, with $n>0$ ($n=0$ is unproblematic). $\underline{C}$ can be written in conjunctive normal form, in which each conjunct can be written as a disjunction of the form ``$x_{n_1}=\alpha_1 \textrm{ or } x_{n_2}=\alpha_2 \textrm{ or } ...$'', with $\alpha_i\in\{0,1\}$. Such disjunctions translate in rules of the form 
``$0\in\{x_i: \alpha_i=0 \}\textrm{ or }1\in\{x_i: \alpha_i=1 \}$'', 
that is ``$\{x_i: \alpha_i=0 \}\subseteq\{1\} \Rightarrow \{x_i: \alpha_i=1 \}\cap\{1\}\not=\emptyset$'', 
that is statements of the same form as those used in a regularity rule when $\mathcal{D}_p=\mathcal{D}_c=\{1\}$. Put simply, the conjunctive normal form of $\underline{C}$ provides the conclusion regularity rule for $C$, following the format of Theorem~\ref{th:truthconstraintforregconn}. The premise regularity rule is provided by the conjunctive normal form for the negation of $C$.
\end{myproof}

Furthermore, a somewhat converse result also holds:
\begin{theorem}\label{th:bivconnarereg}
Every regular connective in $N$-valued logics with a weakly bivalent truth-function (one which takes bivalent values for bivalent inputs) shares a regularity rule with a connective from classical logic.
\end{theorem}

\begin{myproof}
Suppose that $C$ is a regular, bivalent connective. It is clear that the rule of Theorem~\ref{th:truthconstraintforregconn} holds not only for $\underline{C}$, but also for $\widetilde{\underline{C}}$, the restriction of $\underline{C}$ to inputs in $\{0,1\}$. This last truth-function thus satisfies the same regularity rule as the original one, which is thus one of a classical truth-function: the output values are in $\{0,1\}$, by hypothesis, and the question of whether a value belongs to some set of designated values or not is the same question as whether the value is $1$ or not ($\{1\}$ is the set of designated values for classical logic).
\end{myproof}

Finally, we show a result of expressive completeness: a G-conditional is sufficient to deliver all regular rules from classical connectives.

\begin{theorem}\label{th:conditionalissufficient}
Assume that a logic is constant expressive and has a G-conditional. Then for any classical connective $C$, the logic has a regular connective $C^+$ that shares a regularity rule with this connective $C$.
\end{theorem}

\begin{myproof}
The logic has a G-conditional and is constant expressive, so it also has a G-conjunction and a G-disjunction (Theorem~\ref{th:combinationsforconditionals}). 
Let us write the conclusion regularity rule of a classical $n$-ary connective $C$ as usual with a set $\mathcal{B}^c$ of pairs $(B_p,B_c)$ of sets of indices. Consider the connective $C^+$ with a truth-function defined as:
$\underline{C^+}(x_1, ..., x_n)=\bigwedge_{(B_p,B_c)\in\mathcal{B}^c} [ \bigwedge\{x_i:i\in B_p\} \to \bigvee\{x_i:i\in B_c\}]$ (here using the relevant truth-functions for conjunctions, disjunctions and conditionals, including for $\bigwedge$ and $\bigvee$, for which order does not matter and can be held fixed). This connective surely satisfies the conclusion regularity rule; for all $i$ an index of one mixed truth-relation in the intersective mixed truth-relation:
$\underline{C^+}(x_1, ..., x_n)\in \mathcal{D}_c^i$
	iff
	$(\bigwedge_{(B_p,B_c)\in\mathcal{B}^c} [ \bigwedge\{x_i:i\in B_p\} \to \bigvee\{x_i:i\in B_c\}] )\in \mathcal{D}_c^i$
	iff
	$\bigwedge_{(B_p,B_c)\in\mathcal{B}^c} ([ \bigwedge\{x_i:i\in B_p\} \to \bigvee\{x_i:i\in B_c\}] \in \mathcal{D}_c^i)$
	(G-conjunction, see~\ref{th:formcond})
	iff
	$\bigwedge_{(B_p,B_c)\in\mathcal{B}^c} ([ \bigwedge\{x_i:i\in B_p\}\in\mathcal{D}_p^i \Rightarrow \bigvee\{x_i:i\in B_c\}\in\mathcal{D_c^i}] $
	(G-conditional, see~\ref{th:formcond})
	iff
	$\bigwedge_{(B_p,B_c)\in\mathcal{B}^c} ([ \{x_i:i\in B_p\}\subset\mathcal{D}_p^i \Rightarrow \{x_i:i\in B_c\}\cap\mathcal{D}_c^i\not=\emptyset] $
	(G-conjunction and G-disjunction, see~\ref{th:formcond}).
This provides the result by Theorem~\ref{th:truthconstraintforregconn}.

As for the premise regularity rule, it will follow because for classical connectives it is constrained by the conclusion regularity rule in the first place. We start with a derivation similar to the one above:
$\underline{C^+}(x_1, ..., x_n)\in \mathcal{D}_p^i$
	iff
	$(\bigwedge_{(B_p,B_c)\in\mathcal{B}^c} [ \bigwedge\{x_i:i\in B_p\} \to \bigvee\{x_i:i\in B_c\}] )\in \mathcal{D}_p^i$
	iff
	$\bigwedge_{(B_p,B_c)\in\mathcal{B}^c} ([ \bigwedge\{x_i:i\in B_p\} \to \bigvee\{x_i:i\in B_c\}] \in \mathcal{D}_p^i)$
	(G-conjunction, see~\ref{th:formcond})
	iff
	$\bigwedge_{(B_p,B_c)\in\mathcal{B}^c} ([ \bigwedge\{x_i:i\in B_p\}\in\mathcal{D}_c^i \Rightarrow \bigvee\{x_i:i\in B_c\}\in\mathcal{D}_p^i] $
	(G-conditional, see~\ref{th:formcond}, the main difference between this derivation and the previous one is that here we apply a different bullet point for this step)
	iff
	$\bigwedge_{(B_p,B_c)\in\mathcal{B}^c} ([ \forall i\in B_p: x_i\in\mathcal{D}_c^i] \Rightarrow [\exists i\in B_c: x_i\in\mathcal{D}_p^i])$
	(G-conjunction and G-disjunction, see~\ref{th:formcond}).
Although it uses the pairs $(B_p,B_c)$ from the $\mathcal{B}^c$ in the conclusion regularity rule, this last statement can be translated in the premise regularity rule: for classical connectives, the premise and conclusion regularity rules are dual of one another: whenever $C(x_1,..., x_n)\not\in\{1\}$ (that is, the premise rule holds), then it cannot be that $C(x_1,..., x_n)\in\{1\}$ (that is, conclusion rule does not hold). Thus, the premise and conclusion regularity rules are such that:
($\bigwedge_{(B_p,B_c)\in\mathcal{B}^c} ([ \forall i\in B_p: \alpha_i] \Rightarrow [\exists i\in B_c: \alpha_i])$ is true
iff
$\bigwedge_{(B_p,B_c)\in\mathcal{B}^c} ([ \forall i\in B_p: \alpha_i] \Rightarrow [\exists i\in B_c: \alpha_i])$ is false),
whenever $\alpha_1, ..., \alpha_n$ is a list of Booleans.
Chaining this at the end of the previous derivation helps translate from the formula with the $\mathcal{B}^c$ of the conclusion regularity rule to a formula with the $\mathcal{B}^p$ of the premise regularity rule, which provides the result by Theorem~\ref{th:truthconstraintforregconn}.

\end{myproof}

\subsection{Regularity: reductions to 4-valued logics and the specific case of $\truthrelation_{\{1,\#_p\},\{1,\#_c\}}$
}
\label{sec:reg=4value}

%\changeEC{Regularity rules are more general than what can be found in classical logic, but we can show that it is not necessary to go beyond 4-valued logic to exhaust all of their power.}
{Some regularity rules are not exemplified by connectives in classical logic, and one may thus like to investigate regular connectives in logic with more truth values. Here we establish that every regularity rule can be exemplified by a connective in 4-valued logic.}
To see this, let us first remind that one motivation for the notion is that regularity allows one to represent a broad class of $N$-valued logics (Definition~\ref{def:regularcompactlogic}) within $4$-valued logic, as stated in Theorem~\ref{thm:reg}.

\begin{definition}[Logic, Compact, Regular]\label{def:regularcompactlogic}
We define a \emph{logic} to be 
{a language equipped with a consequence relation.}
A logic is \emph{compact} if $\Gamma \vdash \Delta$ implies the existence of finite subsets $\Gamma'$ and $\Delta'$ of $\Gamma$ and $\Delta$ respectively, such that $\Gamma'\vdash \Delta'$. A logic is \emph{regular} if it has only regular connectives. \end{definition}

\noindent
Beyond Gentzen's approach, regularity matters because of the following result from \cite{chemla2018suszko}:

\begin{theorem}\label{thm:reg} Every monotonic, compact, and regular logic is semantically representable by means of an at-most 4-valued truth-relational and truth-functional semantics.
\end{theorem}

The fact that at most 4 values are needed to semantically characterize a monotonic consequence relation was also established by \cite{blasio2017inferentially} and by \cite{french2017valuations}, but with no heed paid to truth-functionality. The addition from Theorem \ref{thm:reg} was that when the consequence relation admits only regular connectives, regularity ensures that the reduction to four values is moreover truth-functional, that is compositional on truth values. This is not the case in general, however, and the reduction of a semantics for a monotonic logic to only four values (or three, or two) may otherwise lead to violations of truth-functionality (as originally stressed by \citealt{suszko1977fregean} for the semantic representation of Tarskian logics by means of just two values). This is one aspect in which regularity is an important constraint: it guarantees that a logical connective behaves truth-functionally over a \emph{minimal} semantic representation of its associated consequence relation.

In fact, the previous result can be made more precise, establishing that all regular connectives have a representative in the 4-valued logic equipped with the truth-relation $\truthrelation_{\{1,\#_p\},\{1,\#_c\}}$:

\begin{theorem}
Consider the truth-relation $\truthrelation_{\{1,\#_p\},\{1,\#_c\}}$ in 4-valued logic with truth-values noted $\{1,\#_p,\#_c,0\}$. Every pair of premise/conclusion regularity rules for an $n$-ary connective is represented by a unique truth-functional $n$-ary connective.
\end{theorem}

\begin{myproof}
For every regularity rule, Theorem~\ref{th:truthconstraintforregconn} shows the constraints that need to be satisfied by a truth-function for the connective to obey the regularity rule. These constraints can always be satisfied for the current truth-relation, because belonging or not to $\mathcal{D}_p$ (here $\{1,\#_p\}$) is independent from belonging or not to $\mathcal{D}_c$ (here $\{1,\#_c\}$). (See also \citealp[Theorem~4.15]{chemla2018suszko}). Furthermore, Theorem~\ref{th:limitedmulitplicity} establishes that no two truth-functions can satisfy the same regularity rules in $\truthrelation_{\{1,\#_p\},\{1,\#_c\}}$, because no two truth-values play the same role.

\end{myproof}

\subsection{Summary}

We are now equipped with a definition of the family of consequence relations of interest to us, namely intersective mixed truth-relations, and with a definition of the family of connectives of relevance to us, namely Gentzen-regular connectives. We started with a presentation mostly related to classical logic, and then argued that Gentzen-regularity rules are best understood in some specific $4$-valued logic. In the coming sections, we will focus on four types of connectives: disjunction, conjunction, negation, and conditional. An exhaustive review of the situation will be provided for 3-valued logics (Section~\ref{sec:three}), and for 4-valued logics (Section~\ref{sec:four}). General results for all pure consequence relations, and all order-theoretic relations are presented in Section~\ref{sec:pureandorderresults}. Necessary and sufficient algebraic conditions for the presence of conditionals are eventually obtained for all $N$-valued logics in the final Section~\ref{sec:N}.

\section{Exhaustive search in three-valued logics}\label{sec:three}

\newcommand{\truthtablethreevalues}[9]{
\begin{tabular}[t]{c|ccc}
	 & 1 & \half & 0\\
\hline
1 & #1 & #2 & #3\\
\half & #4 & #5 & #6\\
0 & #7 & #8 & #9\\
\end{tabular}
}

In three-valued logics, there are exactly five intersective mixed relations,
namely: $ss$, $tt$, $st$, $ts$ and $ss \cap tt$. In this section, we consider which G-connectives among negation, conjunction, disjunction, and the conditional, can be defined for them. We distinguish two cases: 
first, we look at maximally or constant expressive languages (with $\top, \#, \bot$ for $1, \half,0$), then we generalize the results to less expressive languages, namely atomic expressive languages (see Definition~\ref{def:constantsandmaximalexpressiveness}).

\subsection{Maximally (or constant) expressive languages}

With the exception of $ss\cap tt$, all intersective mixed relations in 3-valued logic admit  Gentzen-regular connectives. More specifically:

\begin{theorem}
G-disjunctions and G-conjunctions can be found for all 3-valued intersective mixed relations. G-negations and G-conditionals can be found for all intersective mixed relations except $ss \cap tt$, which admits neither. 
\end{theorem}

\begin{myproof}
These results can be proven by hand, but they follow from the companion computer program available at \url{https://arxiv.org/src/1809.01066v1/anc}.
The results are obtained as follows: the computer program lists all possible truth-functions and, for each such truth-function $c$ and each truth-relation in 3-valued logic, it checks whether, for all pairs of sets of truth values $\gamma$, $\delta$, the truth-function satisfies the regularity rules associated with conjunction, disjunction, and the conditional. The result is thus first obtained at the level of truth values, but it can be shown to be similar for whole propositions, and therefore formulae 
(see Theorem~\ref{th:reducetotruthvalues}).
 The exhaustive list of G-connectives for each consequence relation is as follows:

\begin{description}

\item[Conjunctions] All consequence relations share the following (Strong Kleene) conjunction:

\[
\truthtablethreevalues{1}{\half}{0}{\half}{\half}{0}{0}{0}{0}
\]

Compatible with Theorems~\ref{th:systematicmulitplicity} and \ref{th:limitedmulitplicity}, we observe that
this is the unique conjunction for $st$, $ts$ and $ss \cap tt$, while, in addition to these, 
$ss$ allows for all variants of these tables where a $0$ is replaced by $\half$ or the reverse (hence a total of $2^8$ conjunctions for $ss$);
$tt$ allows for all variants of these tables where a $1$ is replaced by $\half$ or the reverse (hence a total of $2^4$ conjunctions for $tt$).

\item[Disjunctions] All consequence relations share the following (Strong Kleene) disjunction:

\[
\truthtablethreevalues{1}{1}{1}{1}{\half}{\half}{1}{\half}{0}
\]

Compatible with Theorems~\ref{th:systematicmulitplicity} and \ref{th:limitedmulitplicity}, we observe that
this is the unique disjunction for $st$, $ts$ and $ss \cap tt$, whereas, in addition to these, 
$ss$ allows for all variants of these tables where a $0$ is replaced (nonuniformly) by $\half$ or the reverse (hence a total of $2^4$ disjunctions for $ss$);
$tt$ allows for all variants of these tables where a $1$ is replaced (nonuniformly) by $\half$ or the reverse (hence a total of $2^8$ disjunctions for $tt$).

\item[Conditionals] The relation $ss \cap tt$ does not have a conditional.
We can provide a direct proof of this rather central and surprising fact. Suppose that $\to$ were a regular conditional for $ss \cap tt$. Then, we can prove that no value would work for $\half\to0$:
(i)~$\half\not\truthrelation_{ss\cap tt} 0$, 
	so $\not\truthrelation_{ss\cap tt} (\half\to0)$, 
	so $\half\to0$ cannot be $1$.
(ii)~$\half, \half\not\truthrelation_{ss\cap tt} 0$, 
	so $\half\not\truthrelation_{ss\cap tt} (\half\to0)$, 
	so $\half\to0$ cannot be $\half$.
(iii)~$\not\truthrelation_{ss\cap tt} 0$, 
	so $\not\truthrelation_{ss\cap tt}\half$,
	so it's not the case that $0 \truthrelation_{ss\cap tt}$ and $\truthrelation_{ss\cap tt}\half$,
	so $(\half\to0) \not\truthrelation_{ss\cap tt}$,
	so $\half\to0$ cannot be $0$.

Moving away from $ss \cap tt$ then, the following tables are conditionals for the other four consequence relations. Compatible with Theorems~\ref{th:systematicmulitplicity} and \ref{th:limitedmulitplicity}, these tables exhaust the possibilities, except for $ss$ and for $tt$, for which the same replacements as above produce other conditionals ($2^2$ for $ss$, and $2^7$ for $tt$).\footnote{See \cite{jeffrey1963indeterminate} for a related (more restricted) result concerning the conditionals internalizing $tt$-validity.}\medskip

\noindent\begin{tabular}{c@{$\quad$}c@{$\quad$}c@{$\quad$}c}
	$ss$: \truthtablethreevalues{1}{0}{0}{1}{1}{1}{1}{1}{1} &
	$tt$: \truthtablethreevalues{1}{1}{0}{1}{1}{0}{1}{1}{1} &
	$st$: \truthtablethreevalues{1}{\half}{0}{1}{\half}{\half}{1}{1}{1} &
	$ts$: \truthtablethreevalues{1}{\half}{0}{1}{\half}{\half}{1}{1}{1} 
\end{tabular}\medskip

\item[Negations] The relation $ss \cap tt$ does not have a negation. 
This follows from the absence of a conditional and the presence of a disjunction (Theorem \ref{th:combinationsforconditionals}). Alternatively, this can be obtained with essentially the same proof as the one above for conditionals: no value can be assigned to $\neg(\half)$. 
For other consequence relations, the admissible negations are given by the last column of the tables for the conditionals, with variability obtained through the relevant replacements for $ss$ and $tt$, again because of Theorems~\ref{th:systematicmulitplicity} and \ref{th:limitedmulitplicity}.
\qedhere
\end{description}

\end{myproof}

We note, however, that the absence of a G-conditional for $ss\cap tt$ fundamentally depends on the multi-conclusion setting and on the premise regularity rule. In other settings a conditional may be well-behaved, although the options may be more restricted than one can think of given the range of 3-valued conditionals that have been studied:

\begin{theorem}
The relation $ss \cap tt$ has a unique conditional satisfying a single-conclusion version of the conclusion-Gentzen regularity rule (i.e.~the deduction theorem): $\forall \Gamma: \Gamma \vdash A \to B$ iff $\Gamma, A\vdash B$. It is defined as follows:\footnote{This conditional corresponds to a three-valued version of the so-called G\"odel implication, see \cite{hajek1998meta}.}
\[
	\truthtablethreevalues{1}{\half}{0}{1}{1}{0}{1}{1}{1}
	\]

However, the relation $ss \cap tt$ has no conditional satisfying a single-conclusion version of the premise-Gentzen regularity rule, which may be expressed as follows (note that this requires us to accept empty conclusion sets): $\forall \Gamma,A,B: \Gamma, A\to B\vdash $ iff ($\Gamma, B \vdash$ and $\Gamma \vdash A$). The previous conditional does, however, satisfy the right-to-left direction of that rule.

\end{theorem}

\begin{myproof} This result follows from the companion program, but we prove pieces of it by hand to illustrate the process:

\begin{itemize}

\item Show that $\Gamma, A \vdash B $ iff $\Gamma \vdash A\to B$. 

For all semantic interpretation, the right-hand side holds 
iff $v(A\to B)=1$ or there is $x\in v(\Gamma)$ lower than $v(A\to B)$.
That is, iff $v(A\to B)=1$ or ($v(A\to B)\not=1$ and there is $x\in v(\Gamma)$ lower than $v(A\to B)$).
That is, iff $v(A)\leq v(B)$ (see truth-table) or there is $x\in v(\Gamma)$ lower than $v(B)$ (in all cases where $v(A\to B)\not=1$, $v(A\to B)=v(B)$, see truth-table).
That is, iff $v(A)\truthrelation v(B)$ or $v(\Gamma)\truthrelation v(B)$.
That is, iff there is $x$ in $v(\Gamma)$ or $v(A)$ lower than $v(B)$, that is, iff the left-hand side holds.

\item Uniqueness: 
For any such a conditional, now in a constant expressive setting, let us denote by $\overline{x}$ the constant proposition with value $x$, for any $x$ a truth value. So, $\top$ is $\overline{1}$, $\bot$ is $\overline{0}$, $\#$ is $\overline{\half}$.
Then, $x\leq y$ iff $x\truthrelation y$ iff $\overline{x}\models\overline{y}$ iff $\models\overline{x}\to\overline{y}$ iff $\truthrelation x\to y$ iff $x\to y=1$. That accounts for the diagonal $1$ in the truth-table, the $1$s below that diagonal, and the absence of $1$ above the diagonal.
Furthermore, $\#, \top \vdash \#$, hence it $\# \vdash \top \to \#$, and it is thus not possible that $1\to \half=0$.
Finally, $\#, \overline{x} \not\vdash \bot$, for $x=1$ or $x=\half$, hence $\# \not\vdash \overline{x}\to \bot$, and therefore
$x\to 0$ cannot be $\half$.

\item Show that if $\Gamma, B\vdash $ and $\Gamma \vdash A$, then $\Gamma, A\to B \vdash $.  
From the hypotheses, it follows that for every semantic interpretation $v$, then 
($0$ belongs to $v(\Gamma)$ or $v(B)=0$) and (there is $x$ in $v(\Gamma)$ such that $x\leq v(A)$ or $v(A)=1$).
If $0$ belongs to $v(\Gamma)$, the result follows easily. So, we may suppose that $0\not\in v(\Gamma)$. It follows that $v(B)=0$ (simplifying the first parenthesis in the condition above), and that $v(A)>0$ (because $\Gamma\vdash A$). Hence, we are in one of the two top cells in the last column of the table, that is $v(A \to B)=0$, and the result follows.
\qedhere

\end{itemize}
\end{myproof}

\subsection{Less expressive languages}

We have answered our question for languages in which all constants are expressible. 
What happens when the language is less expressive?
When a consequence relation admits a Gentzen-regular connective in a maximally expressive setting, the truth-function can be used to define a regular connective in a less expressive setting too. In a sense, respecting the regularity rule in a maximally expressive setting is like respecting it \emph{intensionally}, so the rule will continue to be satisfied whether or not all relevant propositions are actually expressible. 
But what happens if a given connective does not exist in a maximally expressive setting, can we find a connective that would be appropriate in a \emph{less} expressive one? 

To answer this, we shall consider $ss \cap tt$, which is the only intersective mixed relation lacking regular connectives in a fully expressive setting, namely G-negations and G-conditionals. More general results will soon be obtained (Theorem~\ref{th:fullordertheoretic}), but we can here directly show that even in a less expressive setting $ss \cap tt$ admits neither a G-negation nor a G-conditional.

\begin{theorem}\label{th:no negation for sstt}
Even with an atomic expressive semantics, $ss \cap tt$ can never admit a G-negation.
\end{theorem}

\begin{myproof}
Suppose $\neg$ is a G-negation for a given language, with $ss \cap tt$ as the associated consequence relation. Assume that we can find a formula $A$ and a semantic interpretation $v$ such that $v(A)=\half$. $A\vdash A$, because $ss \cap tt$ is reflexive. By Gentzen-regularity, $A ,\neg A\vdash$ and $\vdash \neg A, A$. From the former, it follows from $tt$-validity that $\neg(\half)=0$, and from the latter it follows from $ss$-validity that $\neg(\half)=1$.
\end{myproof}

\begin{theorem}\label{th:no conditional for sstt}
Even with an atomic expressive semantics, $ss \cap tt$ can never admit a G-conditional.
\end{theorem}

\begin{myproof}
Suppose $\to$ is a G-conditional for a given language, with $ss \cap tt$ as the associated consequence relation, choose $p$ and $q$ two distinct atomic propositions. 
\begin{itemize}
\item $p \vdash q, p$ by reflexivity and monotonicity of $ss \cap tt$. Hence, by Gentzen-regularity:  $\vdash p\to q, p$. 
	Hence, for all $v$ st $v(p)\not=1$, $v(p\to q)=1$. From this, it follows that $\half\to0=1$.
\item $p, q \vdash q$ and $p \vdash p, q$, again by reflexivity of $ss \cap tt$; hence, if $\to$ is Gentzen-regular, $p, p\to q \vdash q$.
	From a valuation $v$ in which $v(q)=0$ and $v(p)=\half$, we infer that necessarily, $\half\to0=0$.
\end{itemize}
Assuming truth-functionality, the two constraints above are incompatible.
\end{myproof}

\subsection{Summary}

The upshot is that, in a multi-premise multi-conclusion setting, all mixed consequence relations admit a G-conjunction, a G-disjunction, a G-negation and a G-conditional. The only `strict intersective' mixed consequence relation, $ss\cap tt$, only admits a G-conjunction and a G-disjunction, but neither a G-conditional nor a G-negation.

\section{Exhaustive search in four-valued logics}\label{sec:four}

Regular connectives have a close connection to 4-valued logics: they correspond to all connectives going with a particular truth-relation (see Section~\ref{sec:reg=4value}). However, there are many intersective mixed consequence relations in 4-valued logic, more than in the 3-valued case. In this section, we provide a computer-aided exhaustive search of all of these consequence relations, and of the G-connectives they admit, for languages that are constant expressive.

\newcommand{\truthtablefourvaluescond}[6]{
\begin{tabular}{c|cccc}
	 & 1 & $\#_1$ & $\#_2$ & 0\\
\hline
1 & 1 & $\#_1$ & $\#_2$ & 0\\
$\#_1$ & 1 & #1 & #2 & #3\\
$\#_2$ & 1 & #4 & #5 & #6\\
0 & 1 & 1 & 1 & 1\\
\end{tabular}
}

\begin{theorem}\label{th:4-valued-all}
There are 167 distinct 4-valued intersective mixed consequence relations.
 Among these:
\begin{itemize}
\item 18 admit a G-conditional (and, therefore, a G-negation, a G-disjunction and a G-conjunction). These are the 16 mixed consequence relations as well as the following two for which we provide the truth-tables of their unique conditionals:

\[
\begin{tabular}{cc}
$\truthrelation_{\{1,\#_1\},\{1,\#_1\}} \cap \truthrelation_{\{1,\#_2\},\{1,\#_2\}}$
&
$\truthrelation_{\{1,\#_1\},\{1,\#_2\}} \cap \truthrelation_{\{1,\#_2\},\{1,\#_1\}}$
\\[1ex]
\truthtablefourvaluescond{1}{$\#_2$}{$\#_2$}{$\#_1$}{1}{$\#_1$}
&
\truthtablefourvaluescond{$\#_1$}{1}{$\#_1$}{1}{$\#_2$}{$\#_2$}
\end{tabular}
\]

\item 28 admit G-conjunctions and G-disjunctions, but no G-negation (nor G-conditional).

\item 
	27 admit G-conjunctions, but no G-disjunction (and therefore no G-negation, nor G-conditional);
	and 
	27 others admit G-disjunctions, but no G-conjunction (and therefore no G-negation, nor G-conditional).

\item 
	The remaining 67 admit no G-disjunction, G-conjunction, G-negation nor G-conditional.

\end{itemize}
\end{theorem}

\begin{myproof}
The proof is provided by the companion computer program available at \url{https://arxiv.org/src/1809.01066v1/anc}.
This program explores arguments at the level of truth values, and results are thus obtained for constant expressive semantics (see Theorem~\ref{th:reducetotruthvalues}).
The program does not explore the whole set of truth-tables as in 3-valued logic above, but only a subset in which we can be sure to find regular conditionals, negations, conjunctions and disjunctions, as long as one exists. It does so by restricting attention to connectives of the following form, in which the $\star$ indicates a place for which all choices of truth values were scanned:

\[
\begin{tabular}{cccc}

Negations & Conditionals & Conjunctions & Disjunctions \\
 \begin{tabular}{c|cccc}
	 & Neg \\
\hline
1 & 0 \\
$\#_1$ & $\star$\\
$\#_2$ & $\star$\\
0 & 1\\
\end{tabular}
&
 \truthtablefourvaluescond{$\star$}{$\star$}{$\star$}{$\star$}{$\star$}{$\star$}
 &
 \begin{tabular}{c|cccc}
	 & 1 & $\#_1$ & $\#_2$ & 0\\
\hline
1 		& 1 & $\#_1$ & $\#_2$ & 0\\
$\#_1$ 	& $\#_1$ & $\star$ & $\star$ & 0\\
$\#_2$ 	& $\#_2$ & $\star$ & $\star$ & 0\\
0 		& 0 & 0 & 0 & 0\\
\end{tabular}
&
 \begin{tabular}{c|cccc}
	 & 1 & $\#_1$ & $\#_2$ & 0\\
\hline
1 		& 1 & 1 & 1 & 1\\
$\#_1$ 	& 1 & $\star$ & $\star$ & $\#_1$\\
$\#_2$ 	& 1 & $\star$ & $\star$ & $\#_2$\\
0 		& 1 & $\#_1$ & $\#_2$ & 0\\
\end{tabular}

\end{tabular}
\]\medskip

It can be shown that if there is a relevant Gentzen-regular connective, one can obtain another Gentzen-regular connective of the same kind by making the necessary replacements to obtain the same values as those fixed in the above tables: the reason is simply that the values above make the Gentzen-regularity rules work.

The only thing left to prove then is the uniqueness of the conditionals provided above, but this follows from Theorem~\ref{th:limitedmulitplicity} (one can check that no two truth values play the same role for the relevant truth-relations).
\end{myproof}

In the previous theorem, consequence relations for which the role of the indeterminates would be switched (e.g., $\truthrelation_{\{1,\#_1\},\{1,\#_1\}}$ and $\truthrelation_{\{1,\#_2\},\{1,\#_2\}}$) are considered distinct, even though the difference only hinges on the name given to the indeterminates. By collapsing consequence relations in which those indeterminates play the same role, the numbers decrease slightly:

\begin{theorem}
There are 97 distinct 4-valued intersective mixed relations, if we assimilate those which are identical except for the name given to the indeterminates. Among these:
\begin{itemize}
\item 12 admit a G-conditional (and, therefore, G-negation, G-disjunction and G-conjunction). These include 10 mixed consequence relations as well as the two intersective mixed relations identified above.

\item 15 admit G-conjunctions and G-disjunctions, but no G-negations (nor G-conditionals).

\item 
	16 admit G-conjunctions, but no G-disjunction (and therefore no G-negations, nor G-conditionals); 
	and 27 others admit G-disjunctions, but no G-conjunction (and therefore no G-negations, nor G-conditionals).

\item 
	The remaining 38 admit no G-disjunction, G-conjunction, G-negation nor G-conditional.

\end{itemize}
\end{theorem}

\begin{myproof}
The proof is provided by the companion computer program available at \url{https://arxiv.org/src/1809.01066v1/anc}. The program here is similar to the one used for the proof of the previous Theorem~\ref{th:4-valued-all}, but additionally it compares all consequence relations and prunes those that are similar to another up to permutation of the two indeterminate truth-values.
\end{myproof}

$N$-valued logic therefore exhibits a wide variety of situations. To motivate looking higher up into $N$-valued logics, we note that it is necessary to go at least to $N=5$ to discover a consequence relation with a G-negation but not a G-conditional.

\begin{fact}
The following truth-relation in 5-valued logic admits a G-negation but no G-conditional:
$$
\truthrelation_{\{1,\#_1,\#_2\},\{1,\#_1\}}
\cap
\truthrelation_{\{1,\#_1,\#_3\},\{1,\#_2\}}
\cap
\truthrelation_{\{1,\#_2,\#_3\},\{1,\#_3\}} 
$$

\end{fact}

\begin{myproof}
This truth-relation has been found and demonstrated to have a G-negation and no G-conditional by a computer program. But we can use later results to verify this result.
\begin{description}
\item[There is a G-negation]
We observe that the operator below (provided by the same computer program) is such that
 $\forall x:$ ($x\in\mathcal{D}_p^i$ iff $\neg(x)\in\mathcal{D}_c^i$) and ($x\in\mathcal{D}_c^i$ iff $\neg(x)\in\mathcal{D}_p^i$), and (the proof of) Theorem~\ref{th:disjunctioncondition} shows that this makes it a G-negation.

\begin{center}
\begin{tabular}{r|ccccc}
	$X$ &	1	&	$\#_3$	&	$\#_2$	&	$\#_1$	&	0\\ \hline
$\neg X$ &	0	&	$\#_1$	&	$\#_2$	&	$\#_3$	&	1
\end{tabular}
\end{center}

\item[There is no G-conditional]
The necessary condition DC1 from Definition~\ref{def:disjconjcompatibility} to have a G-conjunction and a G-disjunction is not satisfied: $\{1,\#_1,\#_2\}$ is included in no other set of designated values, yet all of its truth values belong to other (distinct) sets of designated values in the representation.
\qedhere
\end{description}
\end{myproof}

\section{General results for mixed, pure, and order-theoretic relations}\label{sec:pureandorderresults}

In the previous sections, we have presented an exhaustive investigation of truth-relations of all kinds in 3-valued logics and 4-valued logics. Here, we propose a more in-depth look at specific kinds of consequence relations in $N$-valued logics for all $N$s: \emph{mixed} consequence relations (and among them \emph{pure} consequence relations) and \emph{order-theoretic} consequence relations.

\subsection{Mixed and pure consequence relations in $N$-valued logics: classical connectives}\label{sec:pureresults}

\begin{theorem}\label{th:mixedhasclassical}
Every mixed consequence relation admits a G-conditional and all regular connectives from classical logic.
\end{theorem}

\begin{myproof}

Starting from a $\truthrelation_{\mathcal{D}_p,\mathcal{D}_c}$, we will show that it admits all classical regular connectives, by showing that it admits a G-conditional (see Theorem~\ref{th:conditionalissufficient}). 
One could then use a shortcut: $ss$, $st$, $ts$ and $\truthrelation_{\{1,\#_p\},\{1,\#_c\}}$ all admit a G-conditional, and since they exhaust the possibilities of how the two sets of designated values can be in a superset/subset relation, we could reason from them, that is, $\truthrelation_{\mathcal{D}_p,\mathcal{D}_c}$ will behave like one of these, depending on which (if any) of $\mathcal{D}_p$ and $\mathcal{D}_c$ is included in the other.

Another approach is to exhibit a G-conditional, or a pattern to construct one. In the table below, consider that 
$\#_p$ stands for any truth-value that belongs to $\mathcal{D}_p$ but not $\mathcal{D}_c$,
$\#_c$ stands for any truth-value that belongs to $\mathcal{D}_c$ but not $\mathcal{D}_p$,
$1$ stands for any truth-value that belongs to both, and
$0$ for any truth-value that belongs to neither of the sets of designated values. 
The schema is underspecified, e.g., $\#_p$ means any value of the relevant type, and which specific value it is may differ depending on the column, line or output, but the choice does not matter. Crucially, one can always define a proper truth-function from this schema: it never outputs a value that is not $1$ or $0$ (which exist in all logics), or of the same type as one of the input values that produce this output. One can finally verify easily that a conditional constructed from this schema will always respect the G-conditional regularity rules (this is most easily seen from Theorem~\ref{th:truthconstraintforregconn}, or simply because it is a G-conditional for the most complete $\truthrelation_{\{1,\#_p\},\{1,\#_c\}}$).
\[
\begin{tabular}[b]{c|cccc}
	 	& 1 	& $\#_p$ 	& $\#_c$ 	& 0\\ \hline
1		& 1	& $\#_p$	& $\#_c$	& 0\\
$\#_p$	& 1	& $\#_p$	& 1		& $\#_p$\\
$\#_c$	& 1	& 1		& $\#_c$	& $\#_c$\\
0		& 1	& 1		& 1		& 1
\end{tabular}\qedhere\]
\end{myproof}

The next two theorems separate out mixed consequence relations from pure consequence relations, and reveal further details:

\begin{theorem}\label{th:pureisclassical}
Every pure consequence relation admits exactly the same regular connectives as classical logic.
\end{theorem}

\begin{myproof}
Consider a regular connective for a pure consequence relation. It  can be turned into an equivalently regular connective that would be bivalent, by changing all non classical value into $1$ when they belong to the set of designated values, and $0$ otherwise, as these changes comply with the constraints from Theorem~\ref{th:systematicmulitplicity}. Hence, any such connective satisfies a regularity rule also satisfied by a classical connective, by Theorem~\ref{th:bivconnarereg}.
\end{myproof}

\begin{theorem}\label{th:mixedisnotclassical}
Every non-pure mixed consequence relation admits regular connectives with no counterpart in classical logic (as well as all regular connectives from classical logic).
\end{theorem}

\begin{myproof}
Consider a non-pure mixed consequence relation based on $\truthrelation_{\mathcal{D}_p,\mathcal{D}_c}$. If it is non-pure, we can find an element that is in one of the set of designated values and not in the other, call it $\alpha$. The $0$-ary connective with constant value $\alpha$ is a regular connective, that does not correspond to a regular connective in classical logic: let us show this by showing the regularity rules corresponding to the two possibilities for $\overline{\alpha}$, and note that they are not possible regularity rules for $0$-ary connectives in classical logic (see Example~\ref{regrulesclassicalconnectives})
\begin{itemize}
\item $\alpha\in\mathcal{D}_p\setminus\mathcal{D}_c$.
($\Gamma, \overline{\alpha} \vdash \Delta$ iff $\Gamma \vdash \Delta$)
\phantom{\emph{true}}
($\Gamma \vdash \overline{\alpha}, \Delta$ iff $\Gamma \vdash \Delta$)
\item $\alpha\in\mathcal{D}_c\setminus\mathcal{D}_p$.
($\Gamma, \overline{\alpha} \vdash \Delta$ iff \emph{true})
$\phantom{\Gamma\vdash\Delta}$
($\Gamma \vdash \overline{\alpha}, \Delta$ iff \emph{true})
\qedhere
\end{itemize}
\end{myproof}

One may thus wonder whether the converse of Theorem~\ref{th:pureisclassical} would hold, yielding a powerful characterization of all pure consequence relations from the set of regular connectives they admit, as in Conjecture~\ref{conjecture:pureconsclass}. But this conjecture is false and the counter-example has already been mentioned, as one of the intersective mixed consequence relations that admit a conditional in 4-valued logics:

\begin{conjecture}[false]\label{conjecture:pureconsclass}
Every relation which admits exactly the same regular connectives as classical logic is a pure consequence relation.
\end{conjecture}

\begin{counterexample}
The 4-valued logic based on the truth-relation $\truthrelation_{\{1,\#_1\},\{1,\#_2\}}\cap\truthrelation_{\{1,\#_2\},\{1,\#_1\}}$ admits exactly the same regular connectives as classical logic.
\end{counterexample}

\begin{myproof}
First, it admits all the connectives from classical logic, because it admits a G-conditional (see Theorem~\ref{th:conditionalissufficient}). Conversely, assume that $C$ is a regular connective for this logic, we will show that its restriction to classical inputs is classical (what we called `weakly bivalent' earlier) and therefore, by Theorem~\ref{th:bivconnarereg}, follows the same regularity rules. Consider conclusion regularity rules of the usual form, and apply Theorem~\ref{th:truthconstraintforregconn} to $C$ and each of the members in the intersection forming the truth-relation:

	\[\begin{array}{c@{\quad\textrm{ iff }\quad}c}

	\underline{C}(x_1, ..., x_n) \in\{1,\#_2\}
		& 
		\bigwedge\limits_{(B_p,B_c)\in \mathcal{B}^c} 
			{\{x_i: i\in B_p\}\subseteq\{1,\#_1\} \Rightarrow  \{x_i: i\in B_c\}\cap\{1,\#_2\}\not=\emptyset}
			\\

	\underline{C}(x_1, ..., x_n) \in\{1,\#_1\}
		& 
		\bigwedge\limits_{(B_p,B_c)\in \mathcal{B}^c} 
			{\{x_i: i\in B_p\}\subseteq\{1,\#_2\} \Rightarrow  \{x_i: i\in B_c\}\cap\{1,\#_1\}\not=\emptyset}
			\\
	\end{array}\]

\noindent
Consider classical inputs $x_1, ..., x_n$ for now. We can drop the indeterminate values from the right-hand sides. The two right-hand sides from above become equivalent, and we thus obtain that ($\underline{C}(x_1, ..., x_n)\in\{1,\#_1\}$ iff $\underline{C}(x_1, ..., x_n)\in\{1,\#_2\}$), that is, $\underline{C}(x_1, ..., x_n)\in\{1,0\}$. 
\end{myproof}

\subsection{Order-theoretic relations in $N$-valued logics}
\label{sec:ordertheoreticresults}

In 3-valued logics, the order-theoretic relation $ss\cap tt$ stands out as one that does not admit a G-regular negation or conditional. What about order-theoretic relations in general? We start by capitalizing on the previous section about pure consequence relations.

\begin{theorem}
Every regular connective for an order-theoretic relation shares a regularity rule with a connective from classical logic.
\end{theorem}

\begin{myproof}
Theorem~\ref{th:pureisclassical} shows that pure consequence relations only satisfy classical regularity rules,
Theorem~\ref{th:ordertheoric=intersectionpure} shows that order-theoretic relations are intersections of pure consequence relations,
and Corollary~\ref{cor:intersectioncoincideR} shows that for an intersection of truth-relations to admit a connective satisfying some regularity rule, all its members should have a connective satisfying this regularity rule.
\end{myproof}

However, for the restricted set of G-connectives we have been interested in, we can offer a complete generalization for every $N\geq 3$, regardless of whether the semantics is maximally expressive (as soon as it is atomic expressive).

\begin{theorem}\label{th:fullordertheoretic}
Consider $N$-valued logics with $N\geq 3$, and assume at least atomic expressiveness. 
An order-theoretic relation has no G-negation and no G-conditional. 
It has a G-conjunction and a G-disjunction if and only if it is a total order-theoretic relation.
\end{theorem}

\begin{myproof}
All elements of the proof work in the same way. For each connective $C$, we start by using the fact that the order-theoretic relation is reflexive, and state for atomic propositions that $C(p_1,p_2,...)\vdash C(p_1,p_2,...)$. Then we apply first the premise-regularity rule for $C$, and then the conclusion-regularity rule for $C$. The two together constrain the values that $C$ can take.
\begin{itemize}

\item Assume that there is a G-negation $\neg$. By reflexivity, for every atomic formula $p$: $\neg p \vdash \neg p$. The following two consequences are incompatible.
	\begin{itemize}
	\item Premise-regularity rule: $\vdash p, \neg p$. So, $\forall v: \truthrelation v(p), v(\neg p)$. Choosing $v$ such that $v(p)$ is an indeterminate $\#_i$, it follows that $\neg \#_i=1$.
	\item Conclusion-regularity rule: $p, \neg p\vdash $. So, $\forall v: v(p), v(\neg p) \truthrelation $. Choosing the same $v$ as above, it follows that $\neg \#_i=0$.
	\end{itemize}

\item Assume that there is a G-conditional $\to$. (We would not need this part of the proof in a constant expressive setting, for the absence of a G-negation would guarantee the absence of a G-conditional). By reflexivity, for all atomic formulae $p, q$: $(p\to q) \vdash (p\to q)$. The following two consequences are incompatible.
	\begin{itemize}
	\item Premise-regularity rule: $\vdash p, (p\to q)$ and $q \vdash (p\to q)$. 
	So, in particular, $\forall v: \truthrelation v(p), v(p\to q)$. Choosing $v$ such that $v(p)$ is an indeterminate $\#_i$ and $v(q)=0$, it follows that $v(\#_i\to 0)=1$.
	\item Conclusion-regularity rule: $p, (p\to q) \vdash q$. 
	So, $\forall v: v(p), v(p\to q)\truthrelation v(q)$. Choosing the same $v$ as above, it follows that $v(\#_i\to 0)=0$.
	\end{itemize}

\item Assume that there is a G-conjunction $\wedge$. By reflexivity, for all atomic formulae $p, q$: $(p\wedge q) \vdash (p\wedge q)$. The following two consequences are incompatible.
	\begin{itemize}
	\item Premise-regularity rule: $p, q \vdash (p\wedge q)$.
	For any two indeterminates $\#_i$ and $\#_j$, choose $v_{i,j}$ such that $v_{i,j}(p)=\#_i$ and $v_{i,j}(q)=\#_j$. It follows that 
	$v_{i,j}(p), v_{i,j}(p) \truthrelation v_{i,j}(p\wedge q)$.
	That is $\#_i \leq v_{i,j}(p\wedge q)$ or $\#_j \leq v_{i,j}(p\wedge q)$.

	\item Conclusion-regularity rule: $(p\wedge q) \vdash p$ and $(p\wedge q) \vdash q$. 
	For any two indeterminates $\#_i$ and $\#_j$, choose $v_{i,j}$ such that $v_{i,j}(p)=\#_i$ and $v_{i,j}(q)=\#_j$. It follows that 
	$v_{i,j}(p\wedge q) \truthrelation v_{i,j}(p)$ and $v_{i,j}(p\wedge q) \truthrelation v_{i,j}(q)$.
	That is $v_{i,j}(p\wedge q) \leq \#_i$ and $v_{i,j}(p\wedge q) \leq \#_j$.

	\end{itemize}
	Summarizing, there is an $x$ such that ($\#_i\leq x$ or $\#_j\leq x$) and ($x\leq\#_i$ and $x\leq\#_j$). 
	It follows that ($\#_i\leq\#_j$ or $\#_j\leq\#_i$), that is, any two indeterminates are ordered.
	
\item The proof for G-disjunctions is analogous.
\qedhere

\end{itemize}

\end{myproof}

Order-theoretic relations thus show a very stable behavior across the number of truth values and the expressive resources of the language. In short, they do not admit G-conditionals or G-negations, and they rarely admit G-conjunctions or G-disjunctions.

\section{Algebraic characterization of consequence relations admitting G-conditionals}\label{sec:N}

In the previous sections, we have looked at particular (kinds of) consequence relations (3- or 4-valued logics, pure or order-theoretic consequence relations), and asked whether they admitted certain G-connectives. Here, we show that one may start from a G-connective, and seek a general algebraic characterization of the conditions under which a given intersective mixed relation admits this particular G-connective, for all $N$. We proceed with G-disjunctions and G-conjunctions, then examine G-negations to finally get at G-conditionals. Let us stress from the outset that the algebraic characterizations we obtain are simple for disjunctions and conjunctions, but not so for conditionals, even though we thought that conditionals were the best candidate for a connective to closely resemble its consequence relation.

\subsection{Disjunctions}

\begin{definition}[disjunction-compatible]\label{def:disjunctioncompatible}
A list of sets of designated values $\mathcal{D}_1, ..., \mathcal{D}_n$ is called \emph{disjunction-compatible} if 
for all truth values $x$ and $y$, there is a truth value $z$ that belongs to all the sets of designated values to which $x$ belongs, all the sets of designated values to which $y$ belongs, and to no other set of designated values.
\end{definition}

\begin{theorem}\label{th:disjunctioncondition}
The following statements are equivalent:
\begin{itemize}
\item The truth-relational consequence relation $\vdash$ admits a G-disjunction in a constant expressive setting.
\item All the minimal representations of $\truthrelation$ are based on a list of sets of designated values which is disjunction-compatible.
\item One of the representations of $\truthrelation$ is based on a list of sets of designated values which is disjunction-compatible.
\end{itemize}
\end{theorem}

\begin{myproof}
Corollary~\ref{th:formcond} says that a connective $C$ will be a G-disjunction iff its truth-function $f$ is such that for all sets of designated values taking part in a minimal representation, and for all truth-values:
$f(x,y)\in\mathcal{D}$ iff $x\in\mathcal{D}$ or $y\in\mathcal{D}$.
Disjunction compatibility is equivalent to the possibility to define such a function.

\end{myproof}

\subsection{Conjunctions}

\begin{definition}[conjunction-compatible]
A list of sets of designated values $\mathcal{D}_1, ..., \mathcal{D}_n$ is called \emph{conjunction-compatible} if 
for all truth values $x$ and $y$, there is a truth value $z$ that belongs to all the sets of designated values to which both $x$ and $y$ belong, and to no other.
\end{definition}

\begin{theorem}\label{th:conjunctioncondition}
The following statements are equivalent:
\begin{itemize}
\item The truth-relational consequence relation $\vdash$ admits a G-conjunction in a constant expressive setting.
\item All the minimal representations of $\truthrelation$ are based on a list of sets of designated values which is conjunction-compatible.
\item One of the representations of $\truthrelation$ is based on a list of sets of designated values which is conjunction-compatible.
\end{itemize}
\end{theorem}

\begin{myproof}
Corollary~\ref{th:formcond} says that a connective $C$ will be a G-conjunction iff its truth-function $f$ is such that for all sets of designated values taking part in a minimal representation, and for all truth-values:
$f(x,y)\in\mathcal{D}$ iff $x\in\mathcal{D}$ and $y\in\mathcal{D}$.
Disjunction compatibility is equivalent to the possibility to define such a function.
\end{myproof}

\subsection{Conjunction and disjunction}

Interestingly, the conditions on conjunctions and disjunctions only depend on the list of sets of designated values, independently of whether they are repeated, or whether they are premise sets or conclusion sets. We can describe the situation differently, seeing that what is needed is, roughly, that the sets are pairwise distinct, as well as maximally intersecting with one another.

\begin{definition}[disjunction-conjunction compatible]\label{def:disjconjcompatibility}
A list of sets of designated values $\mathcal{D}_1, ..., \mathcal{D}_n$ is called \emph{disjunction-conjunction-compatible} if 
\begin{description}
\item[DC1:] For each set $\mathcal{D}_i$, there is a truth value that belongs only to $\mathcal{D}_i$, or $\mathcal{D}_i$ is included in another, distinct set of designated values.
\item[DC2:] 
%\changeEC{Picking any list of sets of designated values, there is an element that belongs to all of these sets, to sets in which one of them would be fully included, but to no other set of designated values.\nbPE{I agree with reviewer here, I don't really know how to make sense of DC2}}{
For all non-empty sublists of sets of designated values: $\mathcal{D}'_1, ..., \mathcal{D}'_{n'}$, there is an element $x$ that belongs to all of these sets ($\forall i', x\in\mathcal{D}'_{i'}$), as well as to any set $\mathcal{D}_i$ in which one of the $\mathcal{D}'_{i'}$ would be fully included, but to no other of the original sets of designated values $\mathcal{D}_i$s. In other words, $x\in\mathcal{D}_i$ iff $\exists i': \mathcal{D}'_{i'}\subseteq\mathcal{D}_i$.
%}
\end{description}
\end{definition}

\begin{theorem}\label{char:disj/conj}
The following statements are equivalent:
\begin{itemize}
\item The truth-relational consequence relation $\truthrelation$ admits a G-disjunction and a G-conjunction in a constant expressive setting.
\item All minimal representations of $\truthrelation$ are based on a list of sets of designated values which is disjunction-conjunction compatible.
\item One of the representations of $\truthrelation$ is based on a list of sets of designated values which is disjunction-conjunction-compatible.
\end{itemize}
\end{theorem}

\begin{myproof}
Firstly, suppose that the relation admits a G-disjunction and a G-conjunction, then any minimal representation is based on disjunction-conjunction-compatible sets of designated values:
\begin{itemize}
\item The first condition DC1 follows from there being a conjunction (no need for a disjunction). Consider $x_1$, the conjunction of all truth values in $\mathcal{D}_1$ (no matter the order in which the conjunction is taken, which may vary given that the conjunction is a binary operator here):
$x_1=(x_1^1 \wedge (x_1^2 \wedge ( ... \wedge x_1^{n_1})...)$, if $\mathcal{D}_1$ has $n_1$ elements $x_1^1, ..., x_1^{n_1}$.
Clearly, $x_1$ belongs to $\mathcal{D}_1$ (as a conjunction of elements that do).
If $x_1$ also belongs to $\mathcal{D}_2$, then all the $x_1^1, ..., x_1^{n_1}$ also belong to $\mathcal{D}_2$ (for the same reason).
\item As for the second condition DC2, pick a list of sets of designated values, consider the disjunction of the elements made of the conjunctions of all elements in each of these sets. That element satisfies the constraint DC2.
\end{itemize}

Secondly, suppose that the third statement above is true: there is a representation {of $\truthrelation$} based on a disjunction-conjunction compatible list of sets of designated values.
The condition DC2 by itself guarantees disjunction-compatibility / conjunction-compatibility: for any two truth values $x$ and $y$, 
pick the list of sets of designated values to which one or the other or both belong (if any), and you can construct the necessary $z$ for disjunction-compatibility, 
pick the list of sets of designated values to which both belong (if any), and you can construct the necessary $z$ for conjunction-compatibility.
\end{myproof}

\subsection{Negations}

For G-negation, we simply state a set of necessary conditions, which will prove useful later on:

\begin{definition}[negation-necessity]\label{def:negationnecessary}
A representation $\truthrelation_{\mathcal{D}_p^1,\mathcal{D}_c^1} \cap ... \cap \truthrelation_{\mathcal{D}_p^K,\mathcal{D}_c^K}$ satisfies \emph{negation-necessity} if 
\begin{description}
	\item[N1] There is no inclusion between two $\mathcal{D}_p$s or two $\mathcal{D}_c$s.
	\item[N2] If there is a superset/subset relation between, a $\mathcal{D}_p^i$ and a $\mathcal{D}_c^j$, then there is the same superset/subset relation between $\mathcal{D}_p^j$ and $\mathcal{D}_c^i$.
\end{description}
\end{definition}

\begin{theorem}\label{char:negation:necessary}
If an intersective mixed relation admits a G-negation, then all of its minimal representations satisfy negation-necessity.
\end{theorem}

\begin{myproof}
Assume there is a G-negation $\neg$. We will make use of Theorem~\ref{th:truthconstraintforregconn} applied to negation: ($x\in\mathcal{D}_p^i$ iff $\neg x\not\in\mathcal{D}_c^i$) and ($x\in\mathcal{D}_c^i$ iff $\neg x\not\in\mathcal{D}_p^i$).
\begin{description}
	\item[For N1:] Suppose $\mathcal{D}_p^1\subseteq\mathcal{D}_p^2$. Then pick $x\in\mathcal{D}_c^2$. Then $\neg x\not\in\mathcal{D}_p^2$.  Hence, $\neg x\not\in\mathcal{D}_p^1$. Hence, $x\in\mathcal{D}_c^1$. So, $\mathcal{D}_p^1\subseteq\mathcal{D}_p^2$, and $\mathcal{D}_c^2\subseteq\mathcal{D}_c^1$, we can drop $\truthrelation_{\mathcal{D}_p^1,\mathcal{D}_c^1}$ from the representation, contradicting minimal representation.
	\item[For N2:] Suppose, for instance, that $\mathcal{D}_p^1 \subseteq \mathcal{D}_c^2$ (similar proof for $\mathcal{D}_c^1 \subseteq \mathcal{D}_p^2$.). Pick $x\in\mathcal{D}_p^2$.
Then $\neg x\not\in\mathcal{D}_c^2$.
Hence, $\neg x\not\in\mathcal{D}_p^1$. 
Hence, $x\in\mathcal{D}_c^1$ and, abstracting away from the initial choice of $x$, $\mathcal{D}_p^2\subseteq\mathcal{D}_c^1$.
\qedhere
\end{description}
\end{myproof}

\subsection{Conditionals}

Seeking necessary and sufficient conditions for the existence of a G-conditional, we could start from the form a G-conditional ought to have, as described in Theorem~\ref{th:formcond}. Instead, we capitalize on the fact that allowing for a G-conditional is equivalent to allowing for any of the following set of G-connectives (see Theorem~\ref{th:combinationsforconditionals}): 
(a G-negation and a G-disjunction),
(a G-negation and a G-conjunction), 
(a G-negation, a G-disjunction and a G-conjunction).
As a consequence, we can combine the previous exploration of conditions of existence for a G-disjunction and a G-conjunction (Corollary~\ref{char:disj/conj}) and of conditions for a G-negation (Theorem~\ref{char:negation:necessary}), and obtain necessary and sufficient conditions for the existence of a G-conditional:

\begin{theorem}\label{theorem:fullalgcond}
An intersective mixed relation admits a G-conditional iff all of its minimal representations satisfy disjunction-conjunction-compatibility (Definition~\ref{def:disjconjcompatibility}) and negation-necessity (Definition~\ref{def:negationnecessary}). We repeat these conditions here:
\begin{description}
	\item[DC1:] For each set $\mathcal{D}_i$, there is a truth value that belongs only to $\mathcal{D}_i$, or $\mathcal{D}_i$ is included in another, distinct set of designated values.
	\item[DC2:] %\changeEC{Picking any list of sets of designated values, there is an element that belongs to these sets, to sets in which one of them would be fully included, and to no other set.}{
For all non-empty sublists of sets of designated values: $\mathcal{D}'_1, ..., \mathcal{D}'_{n'}$, there is an element $x$ that belongs to all of these sets ($\forall i', x\in\mathcal{D}'_{i'}$), as well as to any set $\mathcal{D}_i$ in which one of the $\mathcal{D}'_{i'}$ would be fully included, but to no other of the original sets of designated values $\mathcal{D}_i$s. In other words, $x\in\mathcal{D}_i$ iff $\exists i': \mathcal{D}'_{i'}\subseteq\mathcal{D}_i$.%}
	\item[N1] There is no inclusion between two $\mathcal{D}_p$s or two $\mathcal{D}_c$s.
	\item[N2] If there is a superset/subset relation between, a $\mathcal{D}_p^i$ and a $\mathcal{D}_c^j$, then there is the same superset/subset relation between $\mathcal{D}_p^j$ and $\mathcal{D}_c^i$.
\end{description}
\end{theorem}

\begin{myproof}
These conditions are necessary and sufficient for the existence of a G-disjunction and a G-conjunction by Theorem~\ref{char:disj/conj}, and necessary for the existence of a G-negation by Theorem~\ref{char:negation:necessary}. To prove that they are plainly necessary and sufficient for the existence of a G-conditional, we thus need to prove that they are sufficient for the existence of a G-negation. 
Assuming that the conditions are met, we will show that we can meet the requirements imposed by Theorem~\ref{char:negation:necessary} for a G-negation, namely that for all $x$, there is $y$ such that 
($y \in \mathcal{D}_c^i$ iff $x \not\in \mathcal{D}_p^i$)
and
($y \in \mathcal{D}_p^i$ iff $x \not\in \mathcal{D}_c^i$).

Consider a minimal representation of the intersective mixed relation and let $x$ be a truth value. We would like to find $y$ satisfying the conditions above.
The desiderata then is that $y$ belongs to a given list of sets of designated values $\mathcal{D}_1, ..., \mathcal{D}_n$. Consider $y=\alpha_1 \vee ... \vee \alpha_n$ (the order for the disjunction does not matter), with $\alpha_k$ being the conjunction of all elements in $\mathcal{D}_k$. Surely, $y$ is in all the relevant sets: each $\alpha_k$ is in $\mathcal{D}_k$ (as a conjunction of elements that are), and their disjunction $y$ is therefore in all $\mathcal{D}_k$s (because at least one element in the disjunction is). 

Now, conversely, we should prove that $y$ does not belong to an unwanted set of designated values. Suppose $y$ belongs to some set of designated values, let us say a premise set $\mathcal{D}_p^*$ (for a conclusion-set, the proof would proceed in the same manner). Then, there is an $\alpha_k$ such that $\alpha_k\in\mathcal{D}_p^*$ (because of the disjunctive nature of $y$). Hence, $\mathcal{D}_k\subseteq\mathcal{D}_p^*$ (because of the conjunctive nature of $\alpha_k$). One option is that $\mathcal{D}_k$ is $\mathcal{D}_p^*$, which is fine for our needs for $y$. If not, by N1 from Definition \ref{def:negationnecessary}, it follows that $\mathcal{D}_k$ is a conclusion-set in the representation, call it $\mathcal{D}_c^1$. So we have $\mathcal{D}_c^1\subseteq\mathcal{D}_p^*$, and by N2 from Definition \ref{def:negationnecessary}, we also have $\mathcal{D}_c^*\subseteq\mathcal{D}_p^1$.
Now, the reason why $\mathcal{D}_c^1$, aka $\mathcal{D}_k$ was in the list of desiderata was that $x\not\in\mathcal{D}_p^1$, from which it now follows that $x\not\in\mathcal{D}_c^*$, and it is therefore fine that $y$ belongs to $\mathcal{D}_p^*$. That is, we have proved that any $\mathcal{D}^*$ to which $y$ belongs is actually one in the list of the $\mathcal{D}_1,..., \mathcal{D}_n$.
\end{myproof}

Theorem~\ref{theorem:fullalgcond} offers a complete, if not simple, algebraic characterization of the truth-relations allowing for G-conditionals. Let us say again that, with constant expressiveness, these are also conditions under which a truth-relation admits regular connectives corresponding to all classical truth-functions.
We have seen that all mixed consequence relations do admit all classical connectives; some (properly) intersective mixed relations do as well, but not all, and here we can tease apart those that do from those that don't. 

If Gentzen-regularity can be motivated as a criterion of logicality for a set of connectives, by extension this result may be used in favor of the logicality of mixed consequence relations, as well as all others which may admit a conditional. In previous work (\citealt{chemla2017charac}), we had %\changePE{set the order-theoretic relation on a par with pure consequence relations}
{treated order-theoretic relations as equally worthy of being called logical as pure consequence relations}, but based on criteria that did not take into consideration the interaction of logical consequence relations with sentential connectives, namely on the fact that {order-theoretic relations are} reflexive and transitive. In light of our results, the lack of an internalizing conditional for order-theoretic relations may be used to single those out: they satisfy classical structural properties such as reflexivity and transitivity, but the non-classical truth values may be used to allow or prevent the existence of certain connectives, chief among which are conditionals.

\section{Conclusions and open issues}

We have presented general results concerning the definability of Gentzen-regular connectives for intersective mixed consequence relations. In the 3-valued and 4-valued case, moreover, we can enumerate which intersective relations admit Gentzen-regular conjunctions, disjunctions, negations, and conditionals, and exhibit the truth-tables for those operators. For classical logic and some $4$-valued logics, and in $N$-valued logics for mixed and pure consequence relations as well as for order-theoretic consequence relations, we have offered general results that inform us either about these logics, or about the Gentzen-regular connectives.

The value of those results is not just combinatorial. They are a step toward a better understanding of the nature of logic and the way in which consequence relations and logical constants interact. To buttress this point, we hereby mention further directions in which the present framework, its results and its methods, may be fruitfully extended:

\begin{itemize}

\item What happens with other Gentzen-regular connectives? We have explored in depth the situation for four types of binary connectives {but there are in principle much more of them (a subset of the $4^{4 \times 4}$ of all truth-functional connectives in the 4-valued logic described in Section \ref{sec:reg=4value}).} If all of these connectives turn out to be equally important, then the current exploration should be pursued. Alternatively, one may find reasons to focus on a subclass of these connectives: the set we looked at is familiar, albeit arbitrary. One may investigate the whole set of connectives in classical logic for instance, or some other subset of connectives defined by putting algebraic constraints on what may count as a regularity rule (in fact, classical connectives are so constrained because premise and conclusion regularity rules are essentially the conjunctive normal forms of two formulae that are negations from one another, see the proof of Theorem~\ref{th:classicalallregular}). Our hope is that if some subclass of Gentzen regular connectives is found to be important, it should be possible to reproduce our analysis and study this subclass systematically, possibly in a computer-aided fashion as we have proposed.

\item Similarly, we have focussed attention on determining what class of logics admit a particular G-connective, e.g., a conjunction or a conditional. One can ask a broader question: what class of logics admit a particular \emph{set} of G-connectives? For this reason, in the (false) Conjecture~\ref{conjecture:pureconsclass} we asked whether pure consequence relations could be characterized as the relations which admit classical regular connectives and no other regular connectives.
From the $4^{4^n}$ possible G-connectives of arity $n$, not all combinations are possible (e.g., when there is a G-conditional, there is a G-negation). 
One can raise a functional completeness problem: which subset is possible (we know they would have to be clones, see \citealp{KERKHOFF2014107})? More relevant to our current enterprise, what rank is needed to admit exactly a given subset: to obtain all regular connectives, a 4-valued logic is sufficient, what about various proper subsets of all the Gentzen regular connectives? We have seen that the first logic with a G-negation but no G-conditional is found with $N=5$; this pushes the limit of previous results, which tried to reduce the minimal number of truth-values to just $4$ (\citealt{chemla2018suszko}).

\item Conversely, what happens when Gentzen-regularity is relaxed? The methods here employed can serve to refine the conditions put on logical connectives. While Gentzen-regularity is a natural criterion on the interplay between a connective and a consequence relation, for some applications we may want to investigate which consequence relations admit conditional operators satisfying weaker constraints (for example: only the full deduction theorem), or specific constraints (not necessarily implied by Gentzen-regularity).

\item What happens in settings with limited expressiveness? In this paper, we have assumed logics that are expressively complete. We have shown that when expressive completeness is relaxed the admissible Gentzen-regular connectives stay the same in 3-valued logics and, in fact, for all $N$-valued order theoretic relations. However, it remains an open question to discover how different the situation will be if we consider \emph{all} intersective mixed relations.

\item What happens if we move to settings that are single-conclusion, or single-premise? We have seen that, with three truth values already, things change in a single-conclusion setting, in particular $ss\cap tt$ admits a conditional satisfying at least part of the full Gentzen-regularity condition. Since multi-conclusion inference is sometimes viewed with suspicion (\citealt{steinberger2011conclusions}), it is natural to explore the issues we have covered in that setting. {Note, however, that this may require adjusting the regularity rules, which often rely on the possibility to add up conclusions.}

\item 
Finally, while we have started from consequence relations and looked for associated G-connectives, one may consider the reverse problem: given a Gentzen-regular sentential operator, which consequence relations admit it (see Theorem~\ref{th:fromcondtoLC} for a statement of unicity but only under some conditions)?
The question is of importance to determine how tightly a conditional operator is associated to a specific consequence relation.

\end{itemize}

\bibliographystyle{apalike}
\bibliography{biblio}

\end{document}